\documentclass[12pt]{article}

\usepackage{amssymb}
\usepackage{indentfirst}
\usepackage{amssymb,epsfig,amsmath,accents,setspace}
\usepackage{theorem}
\usepackage[usenames,dvipsnames]{color}
\usepackage[active]{srcltx}
\usepackage{bm}
\usepackage{multirow}
\usepackage{lscape}
\usepackage{epstopdf}
\usepackage[justification=centering]{caption}

\usepackage{tikz}

\usepackage{url}
\usepackage[breaklinks]{hyperref}
\hypersetup{
    colorlinks=true,       
    linkcolor=NavyBlue,          
    citecolor=blue,        
    filecolor=magenta,      
    urlcolor=cyan           
}

\usepackage[sort&compress]{natbib}
\bibpunct{(}{)}{;}{a}{,}{,}

\newtheorem{assumption}{Assumption}

\theoremstyle{plain} { \theorembodyfont{\rmfamily}

}

\newtheorem{lemma}{Lemma}
\newtheorem{corollary}{Corollary}
\newtheorem{theorem}{Theorem}

\newcommand{\cF}{{\cal F}}

\newcommand{\bbR}{{\mathbb R}}

\newcommand{\bbX}{{\mathbb X}}

\newcommand{\prob}{\mathbb{P}}

\newcommand{\expect}{\mathbb{E}}

\newcommand{\tran}{^\top}

\def\halmos{\mbox{\quad$\square$}}

\newenvironment{pfof}[1]{\vspace{1ex}\noindent{\em Proof of #1}\hspace{0.5em}}
    {\vspace{1ex}}

\usepackage{setspace}
 \numberwithin{equation}{section}

 \numberwithin{equation}{section}

\begin{document}
\title{Some Distributional Properties of Linear Stochastic Differential Equations
\thanks{The authors acknowledge financial support from the General Research Fund of the Research Grants Council of Hong Kong SAR (Project No. 14200917).}
}
\author{Xue Dong He\thanks{Corresponding Author. Room 505, William M.W. Mong Engineering Building, Department of Systems Engineering and Engineering Management, The Chinese University of Hong Kong, Shatin, N.T., Hong Kong, Telphone: +852-39438336, Email: xdhe@se.cuhk.edu.hk.} and Zhaoli Jiang\thanks{Department of Systems Engineering and Engineering Management, The Chinese University of Hong Kong, Shatin, N.T., Hong Kong, Email: zljiang@se.cuhk.edu.hk.}}
\maketitle
\begin{abstract}
In this paper, we prove a sufficient and necessary condition for the transition probability distribution of a general, time-inhomogeneous linear SDE to possess a density function and study the differentiability of the density function and the transition quantile function of the SDE. Moreover, we completely characterize the support of the marginal distribution of this SDE.

\medskip



\end{abstract}

\section{Introduction}
Stochastic differential equations (SDEs) have wide applications in various fields and linear SDEs are one of the most important classes of SDEs. Examples of applications of linear SDEs include, but are not limited to, the wealth process associated with an affine trading strategy in the Black–Scholes market \citep{HeEtal2019:MedianMaximization}, state dynamics in stochastic linear-quadratic control \citep{YongJZhouXY:99sc}, and physical systems subject to linear fluctuations \citep{risken1985fokker}.

One question that interests researchers in the study of an SDE is the distribution of the SDE at given future time, namely the transition probability distribution of the SDE. For some special linear SDEs, this distribution can be obtained in closed form, examples being the Ornstein-Uhlenbeck process and exponential functionals of Brownian motion with a drift \citep{yor2001exponential}. For a general linear SDE, however, there is no closed form for the transition probability distribution.


A general method to study the transition probability distribution of an SDE is the Kolmogorov forward and backward equations. To apply this method, it is crucial to assume certain non-degeneracy conditions so that the transition probability distribution has a density function. As reviewed later on, those conditions do not hold for a general linear SDE, so the method of Kolmogorov equations is not always applicable in this case.

Another question that interests researchers is the support of an SDE and the set of states that the SDE can reach at certain future time. For example, in the study of time-inconsistent stochastic control problems, the set of reachable states of an SDE is an important constituent of the definition of equilibrium strategies in \citet{HeJiang2019:OnEquilibriumStrategies}.


In the present paper, we consider a general time-inhomogeneous linear SDE whose coefficients are piece-wise continuous with respect to the time variable. The state space of this SDE is one-dimensional, but the Brownian motion that drives the SDE is multi-dimensional. We prove a sufficient and necessary condition for the transition probability distribution of the SDE to possess a density function and study the joint continuity and differentiability of the density function with respect to the initial data and the terminal state. We then establish a result of regularity of the transition quantile function. On the other hand, by generalizing Strock-Varadhan's support theorem to handle the linear SDE in our setting and solving an associated deterministic optimal control problem, we derive in closed form the set of reachable states of the SDE at each time in the future. Finally, we apply our result to an SDE that arises from portfolio selection and derive additional distributional properties of the SDE by exploiting certain special structures of the SDE.

The remainder of the paper is organized as follows: In Section \ref{se:LiteratureReview} we review the literature and compare our results with those in the literature. In Section \ref{linear SDE} we present our main results. In Section \ref{portfolio selection} we apply our results to an SDE arising from portfolio selection. Two technical lemmas and all proofs are presented in the Appendix.

\section{Literature Review}\label{se:LiteratureReview}
\subsection{Literature on Regularity of Transition Probability}


\citet{hormander1967hypoelliptic} studies the issue of when a second order differential operator with smooth coefficients on a manifold is so-called hypoelliptic. He proposes a sufficient condition, usually referred to as H\"ormander's hypothesis; see theorem 1.1 in \citet{hormander1967hypoelliptic}, and (5.6), (5.7) in \citet{williams1981begin} for the detailed form of H\"ormander's hypothesis. \citet{ichihara1974classification} apply the result obtained by \citet{hormander1967hypoelliptic} to study transition probability of a diffusion process, e.g., an SDE, and to this end, one needs to include the differential operator in the time variable, and this leads to a form of H\"ormander's hypothesis applicable to the probability theory; see (5.8) and the discussion following (5.7) in \citet{williams1981begin} and also see Condition (H) on page 128 of \citet{NualartD:06mc}. This form of H\"ormander's hypothesis can be applied to time-homogeneous SDEs only. To deal with time-inhomogeneous SDEs, \citet{ichihara1974classification} consider another form of H\"ormander's hypothesis; see the paragraph preceding (5.8) of \citet{williams1981begin} and equation (1.6) of \citet{cattiaux2002hypoelliptic} for details. Following \citet{hopfner2017strongly}, we name this form weak H\"ormander's hypothesis, and for time-homogeneous SDEs, this form is equivalent to the original form of H\"ormander's hypothesis in (5.8) of \citet{williams1981begin}. 


The weak H\"ormander's hypothesis requires the drift and diffusion coefficients of the SDE to be smooth in the time variable. To weaken this requirement, \citet{chaleyat1984hypoellipticity} propose the so-called restricted H\"ormander's hypothesis; see (1.7) of \citet{cattiaux2002hypoelliptic} for the detailed form of this hypothesis. When drift and diffusion coefficients are indeed smooth, the restricted H\"ormander's hypothesis implies the weak H\"ormander's hypothesis. \citet{kusuoka1984applications} propose a strong ellipticity condition, which implies the restricted H\"ormander's hypothesis.

\citet{florchinger1990malliavin} attempt to prove that H\"ormander's hypothesis, in the form of Condition (H) on page 128 of \citet{NualartD:06mc} and originally applicable to time-homogeneous SDEs only, can be applied to time-inhomogeneous SDEs as well. \citet{cattiaux2002hypoelliptic}, however, point out that there is a flaw in the proof by \citet{florchinger1990malliavin}. It is pointed out by \citet{cattiaux2002hypoelliptic} that the dependence on the time variable in diffusions poses nontrivial difficulties compared to the case of time-homogeneous diffusions. On the other hand, \citet{derridj1971probleme} prove that the weak H\"ormander's  hypothesis is almost necessary when the SDE has analytic coefficients.

The linear SDE that we study is a time-inhomogeneous one. The results in \citet{ichihara1974classification}, \citet{bally1991connection}, and \citet{hopfner2017strongly}, who assume the weak H\"ormander's hypothesis, cannot apply to this SDE because the the drift and diffussion coefficients in this SDE are not smooth in the time variable and thus do not satisfy the weak H\"ormander's hypothesis; see Section \ref{subse:Model}. On the other hand, \citet{chaleyat1984hypoellipticity} assume restricted H\"ormander's hypothesis globally, namely that the restricted H\"ormander's hypothesis holds for any time and state. The linear SDE that we study, however can be degenerate, namely the diffusion coefficient of the SDE can be zero, at some time and state, so we cannot apply the result in \citet{chaleyat1984hypoellipticity}. For the same reason, we cannot apply the result obtained by \citet{kusuoka1984applications} to our setting either.
 \citet{cattiaux2002hypoelliptic} assume that the restricted H\"ormander's hypothesis holds locally and that the drift and diffusions coefficients are H\"older continuous in the time variable with H\"older index larger than certain threshold; see Theorem 4.3 therein. For the linear SDE under our study, however, the coefficients are  piece-wise continuous in the time variable. 

By assuming a local strong ellipticity condition, \citet{stroock1981malliavin} study the differentiability of the transition density function with respect to the initial data and with respect to the terminal state separately. The joint regularity of the transition probability distribution with respect to the initial data and the terminal state of the process is obtained in the literature only when the weak or restricted H\"ormander's hypothesis holds globally; see for instance Theorem 3' in \citet{ichihara1974classification}, Theorem (38.16) in \citet{RogersWilliams2000:DiffusionsVolumeII}, and the references therein. For the linear SDE under our study, however, the H\"ormander's hypothesis cannot hold globally because the diffusion term of the linear SDE can be degenerate in certain state. Thus, the joint regularity results with respect to the initial data and the terminal state  of the linear SDE obtained in the present paper is new.

\subsection{Literature on the Support of an SDE}

Strock-Varadhan's support theorem is a crucial tool to study the support of the law, in the space of continuous functions, of the solution to an SDE. The first version of the theorem is proved by \citet{stroock1972support}, where the authors assume the diffusion coefficients to be bounded. \citet{gyongy1989stability} study the support theorem for linear SDEs whose diffusion coefficients are not bounded. \citet{gyongy1990approximation} and \citet{gyongy1995approximation} consider general SDEs under an assumption that is weaker than the one assumed in \citet{stroock1972support}. See \citet{ondrejat2018support} for a summary of the relevant literature.

With the help of the support theorem, one can represent the support of the marginal distribution of the SDE, namely the support of the distribution of the SDE at a single time point, by a deterministic optimal control problem. Using the Girsanov transform and the support theorem in \citet{gyongy1990approximation}, \citet{zak2014exponential} prove in Lemma 3.4 therein that a particular three-dimensional SDE has support of its marginal distribution to be the whole space $\bbR^3$. \citet{kunita1976support} show that under the global H\"ormander's  hypothesis, the support of the marginal distribution of a time-homogeneous SDE is the whole space; also see the application of this result in  \citet{meyn1993stability} and  \citet{colonius1999topological}. To our best knowledge, for a general linear SDE, a complete picture of the support of its marginal distribution has not been derived in the literature. In the present paper, we derive such a complete picture by solving the associated deterministic optimal control problem.


\section{Main Results}\label{linear SDE}
\subsection{Model}\label{subse:Model}
We first introduce some notations. For any set $A$ in a metric space, denote its interior as $\mathrm{int}(A)$ and its closure as $\mathrm{cls}(A)$. Fix an interval $[a,b]$. For a metric space $\mathbb{B}$, denote by $\mathfrak{C}([a, b];\mathbb{B})$ the set of continuous functions from $[a,b]$ to $\mathbb{B}$ and denote by $\mathfrak{C}_{\mathrm{pw}}([a, b];\mathbb{B})$ the set of piece-wise continuous functions from $[a,b]$ to $\mathbb{B}$, i.e., the set of functions $\xi$ from $[a,b]$ to $\mathbb{B}$ such that $\xi$ is continuous on $[t_{i-1},t_{i})$ with $\lim_{t\uparrow t_i}\xi(t)$ existing, $i=1,\dots, N$, for certain partition $a=t_0<t_1<\dots, <t_N=b$. Denote
by $\mathfrak{C}^\infty([a, b];\bbR^l)$ the set of infinitely differentiable functions from $[a, b]$ to $\bbR^l$ and by $\mathfrak{H}([a, b];\bbR^l)$ the set of absolutely continuous functions from $[a, b]$ to $\bbR^l$.


Consider a $d$-dimensional standard Brownian motion $W(t):=\big(W_1(t),...,W_d(t)\big)\tran$, $t\ge 0$ that lives on a filtered probability space $(\Omega, \mathcal{F}, \{\mathcal F_t \}_{t\geq 0}, \prob)$ satisfying the usual condition. Fix $T>0$. For a $\mathbb{R}^{l}$-valued diffusion process $X(t),t\ge 0$ on $(\Omega, \mathcal{F}, \{\mathcal F_t \}_{t\geq 0}, \prob)$, denote by ${\cal S}_{X^T}$ the support of $X(t),t\in[0,T]$, conditional of the information at time 0, where $X(t),t\in[0,T]$ is considered to be a random variable taking values in $\mathfrak{C}([0, T];\bbR^l)$. Denote by $\mathbb{S}_{X(t)}$ the support of the distribution of $X(t)$, conditional on the information at time 0.


We are interested in the following linear SDE:
\begin{align}\label{eq:LinearSDE}
\left\{
\begin{array}{l}
 dX(t)=\left(c_0(t) + c_1(t)X(t)\right)dt+\left(c_2(t) + c_3(t)X(t)\right)\tran dW(t), \; t\in [0, T],\\
 X(0)=x_0\in \bbR,
 \end{array}
 \right.
\end{align}
where $c_0,c_1\in \mathfrak{C}_{\mathrm{pw}}([0, T];\mathbb{R})$ and $c_2,c_3\in \mathfrak{C}_{\mathrm{pw}}([0, T];\mathbb{R}^d)$. We want to study two distributional properties of this SDE. First, we are concerned about $\mathbb{S}_{X(t)}$.

Second, we want to study the differentiability of conditional probability distributuion of $X(T)$. To this end, for each $t\in[0,T)$ and $x\in \bbR$, denote by $\tilde X(s;t,x),s\in[t,T]$ the solution to \eqref{eq:LinearSDE} that starts from time $t$ and state $x$, i.e., the solution to the following SDE:
\begin{align}\label{eq:LinearSDE2}
\left\{
\begin{array}{l}
d\tilde X(s;t,x)=\left(c_0(s) + c_1(s)\tilde X(s;t,x)\right)ds+\left(c_2(s) + c_3(s)\tilde X(s;t,x)\right)\tran dW(s), \; s\in [t, T],\\
\tilde X(t;t,x)=x.
\end{array}
\right.
\end{align}
Define
\begin{align}
&F(t,x,y):=\prob(\tilde X(T;t,x)\le y),\; y\in \bbR,\label{eq:TransitCDF}\\
&G(t,x,\alpha):=\sup\{y\in \bbR:F(t,x,y)\le \alpha\},\; \alpha\in(0,1)\label{eq:TransitQuantile}
\end{align}
to be respectively the cumulative distribution function (CDF) and the (right-continuous) quantile function of $\tilde X(T;t,x)$. 

\subsection{A Transformation}\label{subse:Transformation}
We can remove the drift of \eqref{eq:LinearSDE} by an increasing, affine transformation. More precisely, define
\begin{align}\label{eq:lambdaTransform1}
  \lambda_0(t):=\int_0^t c_0(z)e^{\int_z^t c_1(s)ds}dz,\; \lambda_1(t):=  e^{\int_0^t c_1(s)ds},\; t\in [0, T].
\end{align}
Straightforward calculation yields
  \begin{align}
 &X(t) =  \lambda_0(t)+\lambda_1(t) (x_0+  X^*(t)) ,\; t\in [0, T],
  \end{align}
 where
   \begin{align}
    &d X^*(t) =\big( c_2^*(t)+ c_3(t) X^*(t)\big)\tran dW(t),\; t \in [0, T],  \quad X^*(0)=0,\label{MappingBetweenOf SDE FixedInitial}\\
 &c_2^*(t):=c_2(t)/\lambda_1(t)+c_3(t)[x_0+\lambda_0(t)/\lambda_1(t)],\; t \in [0, T].\label{c2 *}
  \end{align}
As a result,
\begin{align}
\mathbb{S}_{X(t)} = \lambda_0(t) + \lambda_1(t) \left(x_0 + \mathbb{S}_{X^*(t)} \right).\label{MappingSupport}
\end{align}

Similarly, defining
\begin{align}\label{eq:lambdaTransform2}
  \tilde \lambda_0(t)=\int_t^T c_0(z)e^{\int_z^T c_1(s)ds}dz,\; \tilde \lambda_1(t)=e^{\int_t^T c_1(s)ds},\; t \in [0, T],
\end{align}
we have
  \begin{align}
\tilde X(T;t,x)& =  \tilde X^*(T;t,\tilde \lambda_0(t)+ \tilde\lambda_1(t) x),
  \end{align}
 where
    \begin{align}
    &d \tilde X^*(s;t,x) =\big(\tilde c_2^*(s)+ c_3(s) \tilde X^*(s;t, x)\big)\tran dW(s), \; s \in [t, T],\quad  \tilde X^*(t;t,x)=x,\label{MappingBetweenSolutionsOf SDE Fixed T}\\
&\tilde c^*_2(s):=c_2(s)\tilde \lambda_1(s)-c_3(s)\tilde \lambda_0(s),\quad s \in [0, T].\label{tilde c2 *}
  \end{align}
Define
\begin{align}
   &F^*(t,x,y):=\prob(\tilde X^*(T;t,x)\le y),\; y\in \bbR,\label{eq:TransitCDFTransform}\\
   & G^*(t,x,\alpha):=\sup\{y\in \bbR:  F^*(t,x,y)\le \alpha\},\; \alpha\in(0,1)\label{eq:TransitQuantileTransform}
\end{align}
to be respectively the cumulative distribution function (CDF) and the (right-continuous) quantile function of $\tilde X^*(T;t,x)$.
Then we have
  \begin{align}\label{CDFandQuantileEquation}
F(t,x,y)=  F^*\left(t,\tilde \lambda_0(t)+\tilde \lambda_1(t)x,y  \right),\quad G(t,x,\alpha)  =   G^*\left(t,\tilde \lambda_0(t)+ \tilde\lambda_1(t)x,\alpha \right).
  \end{align}

The above transformations will be used in the following study of the distributional properties of the SDE \eqref{eq:LinearSDE}.

\subsection{Probability Densities}
We first present a result of when the transition probability distribution of \eqref{eq:LinearSDE} admits a density function.

\begin{theorem}\label{th:Density}
  Suppose $c_0,c_1\in \mathfrak{C}_{\mathrm{pw}}([0, T];\mathbb{R})$ and $c_2,c_3\in \mathfrak{C}_{\mathrm{pw}}([0, T];\mathbb{R}^d)$. Fix  $t\in[0,T)$. 
  \begin{enumerate}
  \item[(i)] Suppose that $\tilde c^*_2(s)=c_3(s)=0,\forall s\in[t,T)$. Then, $F(t,x,y) = \mathbf 1_{y\ge \tilde \lambda_0(t)+\tilde \lambda_1(t)x }$ for any $x,y\in \bbR$.
  \item[(ii)] Suppose that $c_3(s)=0$ for all $s\in [t,T)$ and $\tilde c^*_2(s)\neq 0$ for some $s\in [t,T)$. Then,
        \begin{align*}
          F(t,x,y) = \Phi\left(\frac{y-\tilde \lambda_0(t)-\tilde \lambda_1(t)x}{b_{t}}\right),\quad y\in \bbR, x\in\bbR,
        \end{align*}
        where  $b_{t}:=\sqrt{\int_t^T\|\tilde c^*_2(s)\|^2ds}>0$ and $F(t,x,y)$ is infinitely differentiable in $(x,y)$.
        Moreover,
            \begin{align}\label{eq:SuperPolyVanishing}
      \lim_{|x|\uparrow +\infty}|x|^k|g(t,x,y)|=0,\; y\in \bbR,\quad  \lim_{|y|\uparrow +\infty}|y|^k|g(t,x,y)|=0,\; x\in \bbR
    \end{align}
    holds for any $k\ge 1$ and $g$ to be the partial derivatives of $F(t,x,y)$ with respect to $x$ and $y$ of any order. In addition, $F(t,x,y)$ and its partial derivatives with respect to $x$ and $y$ of any order are bounded in $(x,y)\in \bbR^2$.
  \item[(iii)] Suppose that $\tilde c^*_2(s) + \xi c_3(s)=0,\forall s\in[t,T)$ for some $\xi\in \bbR$ and that $c_3(s)\neq 0$ for some $s\in [t,T)$. Then,
  \begin{align*}
     F(t,x,y) =\begin{cases}
       \Phi\left(\frac{\ln (y-\xi)-\ln(\tilde \lambda_0(t)+\tilde \lambda_1(t)x-\xi)-\bar a_t}{\bar b_{t}}\right)\mathbf 1_{y>\xi},& y\in \bbR, x>\frac{\xi-\tilde\lambda_0(t)}{\tilde\lambda_1(t)},\\
       \mathbf 1_{y\ge \xi}, &  y\in \bbR, x=\frac{\xi-\tilde\lambda_0(t)}{\tilde\lambda_1(t)},\\
        1-\Phi\left(\frac{\ln (\xi-y)-\ln(\xi-\tilde \lambda_0(t)-\tilde \lambda_1(t)x)-\bar a_t}{\bar b_{t}}\right)\mathbf 1_{y<\xi}, & y\in \bbR, x<\frac{\xi-\tilde\lambda_0(t)}{\tilde\lambda_1(t)},
     \end{cases}
  \end{align*}
     where $\bar b_{t}:=\sqrt{\int_t^T\|c_3(s)\|^2ds}>0$ and $\bar a_{t}:=- \frac{1}{2}\bar b_{t}^2$, and $F(t,x,y)$ are infinitely differentiable in $(x,y)\in \bbR^2\backslash\big\{(\frac{\xi-\tilde\lambda_0(t)}{\tilde\lambda_1(t)} ,\xi)\big\}$. Moreover, \eqref{eq:SuperPolyVanishing} holds for any $k\ge 1$ and $g$ to be the partial derivatives of $F(t,x,y)$ with respect to $x$ and $y$ of any order.
  \item[(iv)] Suppose that for any $v=(v_1,v_2)\tran\in \bbR^2$ with $\|v\|=1$, there exists $s\in[t,T)$ such that $v_1 c_3(s)+v_2 \tilde c^*_2(s)\neq 0$. Then, $F(t,x,y)$ is infinitely differentiable in $(x,y)$. 
      Moreover, $F(t,x,y)$ and its partial derivatives with respect to $x$ and $y$ of any order are bounded in $(x,y)\in \bbR^2$.
  \end{enumerate}
\end{theorem}

It is straightforward to see that the four cases in Theorem \ref{th:Density} are mutually exclusive and collectively exhaustive. 
The following corollary, which characterizes when $\tilde X(T;t,x)$ has a density function, is a direct consequence of Theorem \ref{th:Density}.

\begin{corollary}\label{coro:density}
 Suppose $c_0,c_1\in \mathfrak{C}_{\mathrm{pw}}([0, T];\mathbb{R})$ and $c_2,c_3\in \mathfrak{C}_{\mathrm{pw}}([0, T];\mathbb{R}^d)$. Fix $t\in[0,T)$ and $x\in \bbR$.
  \begin{enumerate}
    \item[(i)] Suppose that $\tilde c^*_2(s)+x c_3(s)=0,\forall s\in [t,T)$. Then,
    $\tilde X\left(T;t,\frac{x-\tilde \lambda_0(t)}{\tilde \lambda_1(t)}\right)=\tilde X^*(T;t,x)\equiv x $ and thus does not admit a density function.
    \item[(ii)] Suppose that $\tilde c^*_2(s)+x c_3(s) \neq0$ for some $s\in [t,T)$. Then, 
    $\tilde X\left(T;t,\frac{x-\tilde \lambda_0(t)}{\tilde \lambda_1(t)}\right)=\tilde X^*(T;t,x)$ possesses a smooth density function.
  \end{enumerate}
\end{corollary}

\subsection{Differentiability of CDF and Quantile Functions}
  Define
  \begin{align}\label{eq:T0}
    t^*:=\inf\{t\in[0,T):\tilde c_2^*(s)=c_3(s)=0,\forall s\in[t,T)\}
  \end{align}
  with the convention that $\inf\emptyset = T$ and
  \begin{align}\label{eq:T0prime}
    t_*:=\inf\{t\in[0,t^*):\tilde c_2^*(s)+\xi c_3(s)=0,\forall s\in[t,t^*)\text{ and some }\xi \in \bbR\}
  \end{align}
  with the convention that $\inf\emptyset = t^*$. Then, by the definition of $t^*$, it is easy to see that there exists unique $\xi\in \bbR$ such that $\tilde c_2^*(s)+\xi c_3(s)=0,\forall s\in[t_*,t^*)$.

  For any interval $[a,b)$ and open set $O$ in $\bbR^l$, denote by $\mathfrak{C}^{0,\infty}([a,b)\times O)$ the set of functions $g(t,z)$ from $[a,b)\times O$ to $\bbR$ such that its derivatives with respect to $z$ of any order exist and are continuous in $(t,z)$ on $[a,b)\times O$, denote by $\mathfrak{C}^{1,\infty}([a,b)\times O)$ the set of functions $g(t,z)$ from $[a,b)\times O$ to $\bbR$ such that its first-order derivative with respect to $t$ and its derivatives with respect to $z$ of any order exist and are continuous in $(t,z)$ on $[a,b)\times O$. Denote by $\mathfrak{C}_{\mathrm{pw}}^{1,\infty}([a,b)\times O)$ the set of functions $g(t,z)$ from $[a,b)\times O$ to $\bbR^l$ such that there exists $a=t_0<t_1<\dots t_{N}=b$ with $g\in \mathfrak{C}^{1,\infty}([t_{i-1},t_{i})\times O)$, $i=1,\dots, N$.
Denote by $\mathbb{N}_0=\mathbb{N}\cup \{0\}$, where $\mathbb{N}$ is the set of positive integers.


For any function $g(t,x,y)$ that is differentiable in $t$ and twice differentiable in $x$, we define
\begin{align}
  {\cal A}g(t,x,y) = g_t(t,x,y) + \left(c_0(t) + c_1(t)x\right) g_x(t,x,y)+ \frac{1}{2}\|c_2(t) + c_3(t)x \|^2g_{xx}(t,x,y),
\end{align}
where $g_t$, $g_x$, and $g_{xx}$ denote respectively the first-order derivative of $g$ with respect to $t$, first- and second-order derivatives of $g$ with respect to $x$.

The following theorem provides a complete picture of the transition probability distribution $F(t,x,y)$.

 \begin{theorem}\label{th:GeneralDensity}
  Consider the SDE \eqref{eq:LinearSDE2}, suppose $c_0,c_1\in \mathfrak{C}_{\mathrm{pw}}([0, T];\mathbb{R})$ and $c_2,c_3\in \mathfrak{C}_{\mathrm{pw}}([0, T];\mathbb{R}^d)$, recall $t^*$ and $t_*$ as defined in \eqref{eq:T0} and \eqref{eq:T0prime}, respectively, and recall the unique $\xi\in \bbR$ such that $\tilde c_2^*(s)+\xi c_3(s)=0,\forall s\in[t_*,t^*)$. Consider any partitions $0=t_0<t_1<\dots<t_{m}=t_*<t_{m+1}<\dots<t_n=t^*$ such that $c_0$, $c_1$, $c_2$, and $c_3$ are continuous on $[t_{i-1},t_i)$ with the left-limits at $t_i$ existent, $i=1,\dots, n$. Recall $F$ and $F^*$ as defined in \eqref{eq:TransitCDF} and \eqref{eq:TransitCDFTransform}, respectively, recall $\tilde \lambda_0$ and $\tilde \lambda_1$ as defined in \eqref{eq:lambdaTransform2}, and define $\tilde \xi(t):=\big(\xi-\tilde\lambda_0(t)\big)/\tilde\lambda_1(t),t\in [0,T]$. Then, $\tilde \lambda_0$ and $\tilde \lambda_1$ on each $[ t_{i-1}, t_i)$ can be extended to $\mathfrak{C}^{1}([t_{i-1},t_i])$, $i=1,\dots, n$ and the following hold:
  \begin{enumerate}
    \item[(i)] For each $i=m+1,\dots, n$, $t\in [t_{i-1},t_{i})$, $F^*\in \mathfrak{C}^{1,\infty}\big([t_{i-1},t)\times \big(\bbR^2\backslash\{(\xi,\xi)\}\big) \big)$, $F^*_t\in \mathfrak{C}^{0,\infty}\big([t_{i-1},t)\times \big(\bbR^2\backslash\{(\xi,\xi)\}\big) \big)$, and
                  \begin{align}\label{eq:FKFormulaGM}
          {\cal A}F(t,x,y)=0,\quad \forall (x,y)\in \bbR^2\backslash\left\{\left(\tilde \xi(t),\xi\right)\right\},t\in [t_{i-1},t_{i}).
        \end{align}
        Moreover, for any $\tau \in [t_*,t^*)$, $F_t(t,x,y)$ is bounded in $(x,y)\in \bbR^2\backslash\{(\tilde \xi(t),\xi)\},t\in [t_*,\tau] $  and
       \begin{align}\label{eq:GMDerivativeBound}
         \sup_{t\in [t_*,\tau],(x-\tilde \xi(t))^2+(y-\xi)^2\ge \delta^2}\left|\frac{\partial^{\ell+j+k} F}{\partial t^\ell\partial x^j \partial y^k}(t,x,y)\right|<+\infty.
       \end{align}
for any $\delta>0$, $\ell \in\{0,1\}$ and $j,k\in \mathbb{N}_0$.

    \item[(ii)] For each $i=1,\dots, m$, $F^*$ and thus $F$ belong to $\mathfrak{C}^{1,\infty}\big([t_{i-1},t_{i})\times \bbR^2\big)$ and
       \begin{align}\label{eq:FKFormulaNonGM}
          {\cal A}F(t,x,y)=0,\quad (t,x,y)\in[t_{i-1},t_i)\times \bbR^2.
        \end{align}
        Moreover, for any $\tau' \in [0,t_*)$ and $i,j\in \mathbb{N}_0$, $|\frac{\partial^{i+j}F}{\partial x^i\partial y^j}(t,x,y)|$ is bounded in $(t,x,y)\in [0,\tau']\times \bbR^2$, and $\sup_{t\in[0,\tau'],y\in \bbR}|\frac{\partial^{i+j}F_t}{\partial x^i\partial y^j}(t,x,y)|$ is of polynomial growth in $x$.
    \item[(iii)] $F$ is continuous on $[0,t_*)\times \bbR^2\cup \{(t,x,y)\mid (x,y)\in \bbR^2\backslash\{(\tilde \xi(t),\xi)\}, t\in [t_*,t^*)\}   $ and for any $j,k\in \mathbb{N}_0$, $\frac{\partial^{j+k} F}{\partial x^j \partial y^k}(t,x,y)$ is continuous on $[0,t_*)\times \bbR^2\cup\{t_* \}\times \bbR \times (\bbR \backslash \{ \xi\}) \cup\{(t,x,y)\mid (x,y)\in \bbR^2\backslash\{(\tilde \xi(t),\xi)\}, t\in (t_*,t^*)\} $.  For any $\tau \in [t_*,t^*)$, $\delta>0$, and $j,k\in \mathbb{N}_0$,     \begin{align}\label{eq:GMDerivativeBoundGlobal}
          \sup_{t\in [0,\tau],x\in \bbR, |y-\xi|>\delta}\left|\frac{\partial^{j+k} F}{\partial x^j \partial y^k}(t,x,y)\right|<+\infty,
        \end{align}
        and $\sup_{t\in [0,\tau],|y-\xi|>\delta}\left|\frac{\partial^{j+k} F_t}{\partial x^j \partial y^k}(t,x,y)\right|$ is of polynomial growth in $x$.
    \item[(iv)] For any $(t,x,y)\in [0,t_*)\times \bbR^2\cup \{(t,x,y)\mid (x,y)\in \bbR^2\backslash\{(\tilde \xi(t),\xi)\}, t\in [t_*,t^*)\}  $ with $F(t,x,y)\in (0,1)$, we have $F_x(t,x,y)<0$.
    \item[(v)] For any $x\neq \frac{y-\tilde \lambda_0(t^*)}{\tilde \lambda_1(t^*)}$,
    \begin{align}
      \lim_{t\uparrow t^*,(x',y')\rightarrow (x,y)}F(t,x',y') = F(t^*,x,y) = F^*\left(t^*,\tilde \lambda_0(t^*)+\tilde \lambda_1(t^*)x,y  \right)= \mathbf 1_{y\ge \tilde \lambda_0(t^*) +\tilde \lambda_1(t^*) x}.
    \end{align}
  \end{enumerate}
\end{theorem}

Note from the definition of $t^*$ that for any $t\in [t^*,T)$, $F(t,x,y) = \mathbf 1_{y\ge \tilde \lambda_0(t) +\tilde \lambda_1(t) x}$, so $F$ is not differentiable. For any $t\in [t_*,t^*)$, Theorem \ref{th:GeneralDensity}-(i) shows that $F$ is differential at $t$ and infinitely differentiable in $(x,y)$ except at a singular point. For any $t\in [0,t_*)$, Theorem \ref{th:GeneralDensity}-(ii) shows that $F$ is differential at $t$ and infinitely differentiable in $(x,y)$ in the whole space. Theorem \ref{th:GeneralDensity}-(iii) shows when the $F$ is continuous and differentiable at the boundary $t_*$. Theorem \ref{th:GeneralDensity}-(iv) and -(v) reveal further properties of the transition probability distribution.

The following corollary provides a complete picture of the transition quantile functions.

 \begin{corollary}\label{coro:GeneralQuantiles}
 Suppose the same conditions as assumed in Theorem \ref{th:GeneralDensity} hold and denote
 \begin{align*}
   {\cal D}:=(t,x)\in[0,t_*)\times \bbR \cup \{(t,x)\mid x \neq \tilde \xi(t), t\in [t_*,t^*)\}.
 \end{align*}
  Then, the following are true:
  \begin{enumerate}
  \item[(i)] $G(t,x,\alpha)=\tilde\lambda_0(t)+ \tilde\lambda_1(t)x$ for all $t\in[t^*,T]$, $x\in\bbR$, and $\alpha\in(0,1)$. $G(t,x,\alpha)=\xi$ for all $t\in[t_*,t^*)$, $x=\tilde \xi(t)$, and $\alpha\in(0,1)$.
    \item[(ii)]For any $(t,x)\in {\cal D}$ and $\alpha\in (0,1)$, $G(t,x,\alpha)$ is uniquely determined by
\begin{align*}
  F(t,x,G(t,x,\alpha)) = \alpha
\end{align*}
and continuous in $(t,x,\alpha)$ in ${\cal D}\times (0,1)$. Moreover,
\begin{align*}
  &G(t,x,\alpha)\neq \xi,\quad  \forall x \neq \tilde \xi(t), t\in [t_*,t^*),\alpha\in(0,1),\\
  &F_y(t,x,G(t,x,\alpha))>0,\quad \forall (t,x)\in {\cal D},\alpha\in(0,1),
\end{align*}
and $G(t,x,\alpha)$ is infinitely differentiable in $(x,\alpha)$ with derivatives to be continuous in $(t,x,\alpha)$. In particular, we have
\begin{align*}
  & G_x(t,x,\alpha) = -\frac{F_x(t,x,G(t,x,\alpha))}{F_y(t,x,G(t,x,\alpha))}, \quad (t,x)\in {\cal D},\alpha\in(0,1).
\end{align*}
\item[(iii)] $G\in \mathfrak{C}^{1,\infty}\big([t_{i-1},t_{i})\times \bbR\times (0,1)\big)$ for all $i=1,\dots, m$ and $G\in \mathfrak{C}^{1,\infty}\big(\{(t,x)\mid x \neq \tilde \xi(t), t\in [t_{i-1},t_{i})\} \times (0,1)\big)$ for all $i=m+1,\dots, n$. In particular,
    \begin{align*}
&      G_t(t,x,\alpha) = -\frac{F_t(t,x,G(t,x,\alpha))}{F_y(t,x,G(t,x,\alpha))}, \quad (t,x)\in {\cal D},\alpha\in(0,1).
    \end{align*}
\item[(iv)] For any $x\in \bbR$ and $\alpha\in(0,1)$,
\begin{align*}
  \lim_{t\uparrow t^*,(x',\alpha')\rightarrow (x,\alpha)}G(t,x',\alpha') =G(t^*,x,\alpha)=\tilde\lambda_0(t^*)+ \tilde\lambda_1(t^*)x.
\end{align*}
  \end{enumerate}
\end{corollary}

\subsection{Support of the SDE}
In this section, we focus on the solution to \eqref{eq:LinearSDE} and study the support of its solution. We are also interested in the set of states that are reachable by the SDE at a given time point. More precisely, set of {\em reachable states} of $X$ at time $t$, denoted as $\bbX_t$, is defined as follows:
\begin{align*}
  \bbX_t:=\mathrm{int}(\mathbb{S}_{X(t)})\cup \left\{x\in \partial \mathbb{S}_{X(t)} : \prob\big(X(t)\in B_\delta(x)\cap \partial \mathbb{S}_{X(t)}\big)>0\text{ for all }\delta>0\right\},
\end{align*}
where $B_\delta(x)$ denotes the ball with radius $\delta$ and centered at $x$ and $\partial \mathbb{S}_{X(t)}$ is the boundary of $\mathbb{S}_{X(t)}$. In other words, the $ \bbX_t$ is the union of $\mathrm{int}(\mathbb{S}_{X(t)})$ and the smallest relatively closed subset $A$ of $\partial \mathbb{S}_{X(t)}$ such that $\prob(X(t)\in A) = \prob(X(t)\in \partial \mathbb{S}_{X(t)})$. For an application of the set of reachable sets, see \citet{HeJiang2019:OnEquilibriumStrategies}.

By definition, we have $\prob(X(t) \in \bbX_t) = 1$. Moreover, we have $\mathbb{S}_{X(t)} = \mathrm{cls}(\mathbb{X}_t)$ and $\mathrm{int}(\mathbb{S}_{X(t)}) = \mathrm{int}(\mathbb{X}_t)$. In general, however, $\bbX_t\neq \mathbb{S}_{X(t)}$. For example, if $X(t)$ is a geometric Brownian motion with the starting point $x_0>0$, then $\bbX_t=(0,+\infty)$ and $\mathbb{S}_{X(t)} = [0,+\infty)$.

The following theorem provides a complete characterization of $\bbX_t$.

 \begin{theorem}\label{th:GenralSupport}
  Consider the SDE \eqref{eq:LinearSDE} and suppose that $c_0,c_1\in \mathfrak{C}_{\mathrm{pw}}([0, T];\mathbb{R})$ and $c_2,c_3\in \mathfrak{C}_{\mathrm{pw}}([0, T];\mathbb{R}^d)$. 
  Define $h(s):=-c_2^*(s)\tran \frac{c_3(s)}{\|c_3(s)\|^2}\mathbf{1}_{c_3(s)\neq 0},s\in [0,T]$, $D:=\{s\in [0,T]:c_3(s)\neq 0\}$, and
   \begin{align}
\underline{t}:=&\inf\{s\in [0, T]:c_2^*(s)\neq 0 \}, \\
  \bar{t}:=&\inf\left\{s\in [0, T]: \int_0^s\left\|c_2^*(z)+ c_3(z)h(z)\right\|dz>0  \right\}
  \end{align}
  with the convention $\inf\emptyset:=T$. Then, 
  $\underline{t} \leq \bar{t}$. Moreover, for each $t\in [0,T]$, $\mathbb{S}_{X(t)}$ is an interval with the left- and right- ends denoted as $\underline{x}(t)$ and $\bar x(t)$, respectively, and $\underline{x}$ and $\bar x$ are left-continuous on $[0,T]$. Furthermore, the following hold:
  \begin{enumerate}
    \item[(i)] For each $t\in[0,\underline{t}]$, $\mathbb{X}_t=\{\lambda_0(t)+\lambda_1(t) x_0\}$.
    \item[(ii)] For each $t\in(\bar{t}, T]$, $X(t)$ possesses a density and $\mathbb{X}_t=\bbR$.
    \item[(iii)] Suppose $\underline{t}<\bar t$ and fix any $t\in (\underline{t},\bar{t}]$. Define $\tau_t:=\sup\{s\in [\underline{t}, t) : c_3(s)\neq 0\}$ with $\sup \emptyset:=\underline{t}$. Then, $X(t)$ possesses a density, $\tau_t >\underline{t}$, and the following hold:
         \begin{enumerate}
     \item[(a)]  Suppose there exist $s_1, s_2 \in [\underline{t}, t) $ such that $h(s_1)<0$ and $h(s_2)>0$. Then,  $\mathbb{X}_t=\mathbb{R}$.
     \item[(b)]  Suppose that $h(s)\le 0,\forall s\in [\underline{t},t)$. If $h$ is not decreasing on $[\underline{t},t)\cap D$, i.e., if there exists $s_1,s_2\in [\underline{t},t)\cap D$ with $s_1<s_2$ and $h(s_1)<h(s_2)$, then $\mathbb{X}_t=\mathbb{R}$. Otherwise, $\mathbb{X}_t=\big(\lambda_0(t)+\lambda_1(t) (x_0+  h^*(\tau_t)), +\infty\big)$, where $h^*(\tau_t):= \lim_{[\underline{t},\tau_t)\cap D\ni s\uparrow \tau_t}h(s)$.
     \item[(c)]  Suppose that $h(s)\ge 0,\forall s \in [\underline{t}, t)$. If $h(s)$ is not increasing on $[\underline{t}, t) \cap D$, i.e., if there exists $s_1,s_2\in [\underline{t},t)\cap D$ with $s_1<s_2$ and $h(s_1)>h(s_2)$, then $\mathbb{X}_t=\mathbb{R}$. Otherwise, $\mathbb{X}_t=\big(-\infty, \lambda_0(t)+\lambda_1(t) (x_0+  h^*(\tau_t) )\big)$, where $h^*(\tau_t):= \lim_{[\underline{t},\tau_t)\cap D\ni s\uparrow \tau_t}h(s)$.
   \end{enumerate}
  \end{enumerate}
\end{theorem}

\section{Linear SDE that Arises from Portfolio Selection}\label{portfolio selection}
As an application, we consider following linear SDE
\begin{align}\label{eq:SDEPortSelection}
\left\{
\begin{array}{l}
 dX(t)=\left(\theta_0(t) + \theta_1(t)X(t)\right)\tran b(t) dt  +\left(\theta_0(t) + \theta_1(t)X(t)\right)\tran \sigma(t) dW(t), \quad t\in [0, T],\\
 X(0)=x_0\in \mathbb{R}
 \end{array}
 \right.
\end{align}
that arises from portfolio selection. The parameters satisfy the following assumption:
\begin{assumption}\label{as:PortfolioSelection}
  $\theta_0,\theta_1,b\in \mathfrak{C}_{\mathrm{pw}}([0, T];\mathbb{R}^n)$, $\sigma\in \mathfrak{C}_{\mathrm{pw}}([0, T];\mathbb{R}^{n\times d})$, and $\sigma(t)\sigma(t)\tran$ is positive definite for any $t\in [0,T]$.
\end{assumption}
In the SDE \eqref{eq:SDEPortSelection}, $b$ and $\sigma$ stand for the mean excess return rate and volatility, respectively, of $n$ stocks in a financial market, $\theta_0$ and $\theta_1$ represent an affine investment strategy, and $X$ stands for the discounted wealth process associated with that strategy. In other words, the discounted dollar amount invested in the stocks in an infinitesimally small period at time $t$ is $\theta_0(t)+\theta_1(t)X(t)$. The positive-definiteness of $\sigma(t)\sigma(t)\tran$ is a standard assumption in portfolio selection.

The SDE \eqref{eq:SDEPortSelection} is a special case of \eqref{eq:LinearSDE} with
\begin{align}\label{eq:cPS}
  c_0(t)=  b(t)\tran \theta_0(t),\quad c_1(t) = b(t)\tran \theta_1(t),\quad c_2(t) = \sigma(t)\tran \theta_0(t),\quad c_3(t) = \sigma(t)\tran\theta_1(t).
\end{align}
As a result, $c_2^*$, $\tilde c_2^*$, $\lambda_0$, $\lambda_1$, $\tilde \lambda_0$, and $\tilde \lambda_1$ as defined in Section \ref{subse:Transformation} become
\begin{align}
\lambda_0(t)& = \int_0^t e^{\int_{s}^tb(z)\tran \theta_1(z)dz}b(s)\tran \theta_0(s) ds,\quad \lambda_1(t) = e^{\int_0^tb(s)\tran \theta_1(s)ds},\label{eq:lambdaPS}\\
   c_2^*(t)&=\sigma(t)\tran \left[\theta_0(t)e^{-\int_0^t b(s)\tran \theta_1(s)ds}+ \theta_1(t)\left(x_0+\int_0^t b(s)\tran \theta_0(s) e^{-\int_0^s b(z)\tran \theta_1(z)dz} ds\right)\right],\label{eq:c2starPS}\\
   \tilde \lambda_0(t) &= \int_t^Te^{\int_s^Tb(z)\tran \theta_1(z)dz}b(s)\tran \theta_0(s)ds,\quad \tilde \lambda_1(t) = e^{\int_t^Tb(s)\tran \theta_1(s)ds},\label{eq:tildelambdaPS}\\
   \tilde c_2^*(t) & = \sigma(t)\tran\left[\theta_0(t)e^{\int_t^Tb(s)\tran\theta_1(s)ds} - \theta_1(t)\int_t^Tb(s)\tran \theta_0(s)e^{\int_s^Tb(z)\tran \theta_1(z)dz}ds\right].\label{eq:tildec2starPS}
\end{align}

%

The following two corollaries are obtained by applying our results in Section \ref{linear SDE} to the particular SDE \eqref{eq:SDEPortSelection} and exploiting additional structures of \eqref{eq:SDEPortSelection}.

\begin{corollary}\label{th:PortDensity}
 Consider the SDE \eqref{eq:SDEPortSelection} and suppose Assumption \ref{as:PortfolioSelection} holds. Then, Theorem \ref{th:GeneralDensity} and Corollary \ref{coro:GeneralQuantiles} hold for SDE \eqref{eq:SDEPortSelection} with $c_i,i=0,1,2,3,4$, $\lambda_i,i=0,1$, and $c^*_2$ as given in \eqref{eq:cPS}, \eqref{eq:lambdaPS}, and \eqref{eq:c2starPS}, respectively. Moreover, we have
 \begin{align}
   &t^*=\inf\{t\in[0,T):\theta_0(s)=\theta_1(s)=0,\forall s\in[t,T)\}, \label{eq:tupperstarPS}\\
   &t_*=\inf\{t\in[0,t^*):\text{there exists }\xi \in \bbR\text{ such that }\theta_0(s)+\xi\theta_1(s)=0,\forall s\in[t,t^*)\},\label{eq:tlowerstarPS}
  \end{align}
 and for any $t\in [t_*, T]$, 
 $\tilde \xi(t)=\xi$.
\end{corollary}

\begin{corollary}\label{prop:PortSupport}
Consider the SDE \eqref{eq:SDEPortSelection} and suppose Assumption \ref{as:PortfolioSelection} holds. Then, Theorem \ref{th:GenralSupport} holds for SDE \eqref{eq:SDEPortSelection} with $c_i,i=0,1,2,3,4$, $\lambda_i,i=0,1$, and $c^*_2$ as given in \eqref{eq:cPS}, \eqref{eq:lambdaPS}, and \eqref{eq:c2starPS}, respectively. Moreover, \begin{align}
    \underline{t}&=\inf\{s\in [0, T]: \theta_0(s)+x_0\theta_1(s) \neq 0 \},\label{eq:tunderlinePS} \\
    \bar{t}&=\inf\{s\in [0, T]: \int_0^s\|\theta_0(z)+\tilde h(z) \theta_1(z)\|dz>0  \},\label{eq:tbarPS}
  \end{align}
  where \begin{align}
  &\tilde h(t):=-\frac{ \theta_0(t) \tran \sigma(t)\sigma(t)\tran \theta_1(t)}{\| \sigma(t)\tran \theta_1(t) \|^2}\mathbf{1}_{ \theta_1(t)\neq 0},\quad t\in [0,T],\label{eq:htildefunPS}\\
  &x^*(t):=\int_0^t b(s)\tran \theta_0(s)e^{\int_s^t b(z)\tran \theta_1(z)dz}ds+x_0 e^{\int_0^t b(s)\tran \theta_1(s) ds},\quad t\in [0,T],\label{eq:xstarPS}
\end{align}
$\mathbb{X}_t=\{x_0\}$ and $x^*(t)=x_0$ for any $t\in[0,\underline{t}]$, $\mathbb{X}_t$ is increasing in $t\in[0,T]$, and
\begin{align}\label{eq:hfunPS}
  h(t) = e^{-\int_0^tb(s)\tran \theta_1(s)ds}\left[\tilde h(t) - x^*(t)\mathbf 1_{\theta_1(t)\neq 0}\right],\quad t\in[0,T].
\end{align}
Furthermore, with $\underline{t}<\bar t$ and fixing any $t\in (\underline{t},\bar t]$ with $\theta_1(s)\neq 0,\forall s\in (\underline{t},t)$, the following are true:
         \begin{enumerate}
     \item[(1)]  Suppose  $\tilde h(s)\leq x^*(s)$ for any $s \in  (\underline{t}, t) $. If $\tilde h$ is not decreasing on $(\underline{t}, t) $, then  $\mathbb{X}_t=\mathbb{R}$.
         If $\tilde h$ is decreasing in on $(\underline{t}, t)$, then $\tilde h(s)< x^*(s)$ for any $s \in  (\underline{t}, t) $ and   $\mathbb{X}_t=(\tilde h(t-), +\infty)$, where $\tilde h(t-):=\lim_{s\uparrow t}\tilde h(s)$.
     \item[(2)]  Suppose $\tilde h(s)\geq x^*(s)$ for any $s \in (\underline{t}, t)$. If $\tilde h$ is not increasing on $(\underline{t}, t)$, then $\mathbb{X}_t=\mathbb{R}$. If $\tilde h$ is increasing on $(\underline{t}, t)$, then $\tilde h(s)> x^*(s)$ for any $s \in  (\underline{t}, t) $ and $\mathbb{X}_t=(-\infty, \tilde h(t-))$, where $\tilde h(t-):=\lim_{s\uparrow t}\tilde h(s)$.
   \end{enumerate}
\end{corollary}



\appendix

\section{Two Lemmas}
In this section, we provide two technical lemmas that will be used in the proofs of the main results of the present paper. The proofs of these two lemmas are presented at the end of the section.

\begin{lemma}\label{lemma:Density}
  Suppose $c_2^*,c_3\in \mathfrak{C}_{\mathrm{pw}}([0, T];\mathbb{R}^d)$. Fix $t\in [0,T)$ and define $\tilde Z_1(s;t),s\in[t,T]$ and $\tilde Z_2(s;t),s\in[t,T]$ by
        \begin{align*}
          d\tilde Z_1(s;t)&= \tilde Z_1(s;t) c_3(s) \tran     dW(s), \;s\in[t,T],\quad \tilde Z_1(t;t)=1,\\
d\tilde Z_2(s;t)&=\left(\tilde c^*_2(s)  +c_3(s)     \tilde Z_2(s;t)\right)\tran dW(s), \\
& \qquad s \in [t, T], \quad \tilde Z_2(t;t)=0.
        \end{align*}
Suppose that for any $v=(v_1,v_2)\tran\in \bbR^2$ with $\|v\|=1$, there exists $s\in[t,T)$ such that $v_1 c_3(s)+v_2 \tilde c^*_2(s)\neq 0$. Then, $\big(\tilde Z_1(T;t),\tilde Z_2(T;t)\big)$ admits an infinitely differentiable probability density $g$ on $\bbR^2$ with
\begin{align}\label{eq:JointDenGrowthCon}
\sup_{z\in  \bbR^2}\|z\|^k\left|\frac{\partial^{i+j}g}{\partial z_1^i\partial z_2^j}(z)\right|<+\infty
\end{align}
for any $k\ge 1$, $i\ge 0$, and $j\ge 0$.
\end{lemma}

\begin{lemma}\label{support set multiD}
  Consider the following linear SDE
\begin{align}\label{multi linear SDE}
\left\{
\begin{array}{l}
  dX(t)=\big(\beta_0(t)+\beta_1(t)X(t)\big)dt  +\sum_{i=1}^d \big(\beta_{2,i}(t)+\beta_{3,i}(t)X(t)\big)dW_i(t),\; t \in [0, T], \\
  X(0)=x_0\in \mathbb{R}^l,
  \end{array}
  \right.
  \end{align}
where $\beta_0,\beta_{2,i}\in \mathfrak{C}_{\mathrm{pw}}([0, T]; \bbR^{l})$ and $\beta_1,\beta_{3,i}\in \mathfrak{C}_{\mathrm{pw}}([0, T]; \bbR^{l\times l})$, $i=1,\dots, d$. Denote ${\cal U}:= \{f_{w} \in \mathfrak{C}([0, T]; \bbR^l):   w\in \mathfrak{H}([0,T];\bbR^d)\}$ and $\bar {\cal U}:= \{f_{w} \in \mathfrak{C}([0, T]; \bbR^l):   w\in \mathfrak{C}^\infty([0,T];\bbR^d)\}$, where $f_{w}$ is the solution of following ODE:
  \begin{align}\label{x w multiD}
f_{w}'(t)&=\beta_0(t)-\frac{1}{2}\sum_{i=1}^d \beta_{3,i} \beta_{2,i}(t)+\left[\beta_1(t)-\frac{1}{2}\sum_{i=1}^d \beta_{3,i}(t) \beta_{3,i}(t)\right]f_{w}(t) \nonumber \\
&\quad +\sum_{i=1}^d\big(\beta_{2,i}(t)+\beta_{3,i}(t)f_{w}(t) \big)w'_i(t),
\; t\in [0, T], \quad f_{w}(0)=x_0.
\end{align}
Then, the following are true:
\begin{enumerate}
  \item[(i)] ${\cal S}_{X^T}=\mathrm{cls}({\cal U}) = \mathrm{cls}(\bar {\cal U})$.
  \item[(ii)] For each $t\in [0,T]$, $\mathbb{S}_{X(t)} = \mathrm{cls}(\mathbb{U}_t)= \mathrm{cls}(\bar {\mathbb{U}}_t)$, where $\mathbb{U}_t:=\{f(t):f\in {\cal U}\}$ and $\bar {\mathbb{U}}_t:=\{f(t):f\in \bar {\cal U}\}$.
  \item[(iii)] $\bar {\cal U}$ and $\bar {\mathbb{U}}_t$, $t\in[0,T]$, are connected. Consequently, when $l=1$, $\mathbb{S}_{X(t)}$ is a closed interval for any $t\in[0,T]$.
\end{enumerate}
  \end{lemma}

Lemma \ref{support set multiD}-(i) is an extension of Theorem 4.1 in \citet{gyongy1989stability} by allowing $\beta_0$, $\beta_1$, $\beta_{2,i}$'s and $\beta_{3,i}$'s to be continuous in multiple pieces of $[0,T]$. Lemma \ref{support set multiD}-(ii) and (iii) are direct consequences of Lemma \ref{support set multiD}-(i).

\begin{pfof}{Lemma \ref{lemma:Density}}
  In the following, we prove that $\big(\tilde Z_1(T;t),\tilde Z_2(T;t)\big)$ is a nondegenerate random vector in the sense of Definition 2.1.1 of \citet{NualartD:06mc}. Then, the lemma is just a consequence of Proposition 2.1.5 of \citet{NualartD:06mc}.

We recall some notations in Malliavin calculus: $\langle f,g\rangle$ stands for the inner product of $f$ and $g$ in the Hilbert space of square-integrable functions from $[t,T]$ to $\bbR$. $D$ denotes the Malliavin derivative operator on the space of square-integrable stochastic processes (see p. 25 of \citealt{NualartD:06mc}). Because for each stochastic process $X$, $DX$ is also a stochastic process, so $DX$ can be identified as $D_sX,s\in[t,T]$. $\mathbb{D}^\infty$ denotes the space of certain smooth random variables (see p. 67 of \citealt{NualartD:06mc}). The Malliavin matrix of $\big(\tilde Z_1(T;t),\tilde Z_2(T;t)\big)$ is defined to be
\begin{align*}
\Gamma = \begin{pmatrix}
\langle D \tilde Z_1(T;t), D \tilde Z_1(T;t)\rangle &  \langle D \tilde Z_1(T;t), D \tilde Z_2(T;t)\rangle\\
\langle D \tilde Z_2(T;t), D \tilde Z_1(T;t)\rangle &  \langle D \tilde Z_2(T;t), D \tilde Z_2(T;t)\rangle
\end{pmatrix}.
\end{align*}
Recalling Definition 2.1.1 of \citet{NualartD:06mc}, to prove that $\big(\tilde Z_1(T;t),\tilde Z_2(T;t)\big)$ is nondegenerate, we need to prove that (i) $\tilde Z_i(T;t)\in \mathbb{D}^\infty$, $i=1,2$ and (ii) $\Gamma$ is invertible almost surely and its determinant, denoted as $\mathrm{det}(\Gamma)$, satisfies $\big(\mathrm{det}(\Gamma)\big)^{-1}\in \cap_{p\ge 1}L^p(\Omega)$, where $L^p(\Omega)$ denotes the space of random variables $X$ with $\expect[|X|^p]<+\infty$. Because  $\tilde c^*_2$, and $c_3$ are piece-wise continuous and thus bounded, Theorem 2.2.2 in \citet{NualartD:06mc} implies $\tilde Z_i(T;t)\in \mathbb{D}^\infty$, $i=1,2$. Thus, we only need to prove (ii) in the following.

Theorem 2.2.1 and equation (2.53) in \citet{NualartD:06mc} imply that
\begin{align*}
&D_s \tilde Z_1(T;t)=c_3(s) \tran \tilde Z_1(s;t) \tilde Z_0(s), \quad D_s \tilde Z_2(T;t)= [\tilde c^*_2(s )+c_3(s)\tilde Z_2(s;t) ]  \tran   \tilde Z_0(s)  , \quad s\in [t, T]
\end{align*}
where
\begin{align*}
\tilde Z_0(s) = \exp\left\{\int_s^T -\frac{1}{2}\|  c_3(\tau)\|^2 d\tau +\int_s^T c_3(\tau) \tran dW(\tau) \right\}.
\end{align*}
Then, denoting the component of $\Gamma$ in its $i$-th row and $j$-th column as $\Gamma_{i,j}$, we have
\begin{align*}
\Gamma_{1,1}&=\langle D \tilde Z_1(T;t), D \tilde Z_1(T;t)\rangle=\int_t^T \|c_3(s) \|^2 \tilde Z_1(s;t)^2 \tilde Z_0(s)^2 d s,\\
\Gamma_{2,2}&=\langle D \tilde Z_2(T;t), D \tilde Z_2(T;t)\rangle=\int_t^T \| \tilde c^*_2(s )+c_3(s)\tilde Z_2(s;t) \|^2\tilde Z_0(s)^2 d s,\\
\Gamma_{1,2}&=\Gamma_{2,1}=\langle D \tilde Z_1(T;t), D \tilde Z_2(T;t)\rangle=\int_t^T [\tilde c^*_2(s )+c_3(s)\tilde Z_2(s;t) ] \tran c_3(s) \tilde Z_0(s)^2 d s.
\end{align*}
Then, H\"older's inequality implies that $\mathrm{det}(\Gamma)\ge 0$, so by Lemma 2.3.1 in \citet{NualartD:06mc}, to prove that $\Gamma$ is invertible almost surely and $\big(\mathrm{det}(\Gamma)\big)^{-1}\in \cap_{p\ge 1}L^p(\Omega)$, we only need to prove that for any $p\geq 2$, there exists $\eta(p)>0$, such that for all $\eta \in (0, \eta(p)]$, we have
\begin{align}\label{Malliavin matrix Lemma}
\sup_{\|v\|=1} \prob( v\tran \Gamma v \leq \eta) \leq \eta^p.
\end{align}

For all $v=(v_1, v_2)\tran\in \bbR^2$ with $\|v\|=1$, we have
\begin{align}\label{Malliavin matrix inequality}
v\tran \Gamma v
&=\int_{t}^T   \| c_3(s) (v_1\tilde Z_1(s;t)+v_2 \tilde Z_2(s;t)   )  +v_2\tilde c^*_2(s )   \|^2 \tilde Z_0(s)^2 ds \nonumber \\
&\geq  \left(\inf_{s \in [t, T]}\tilde Z_0(s)^2\right)\int_{t}^T   \| c_3(s) (v_1\tilde Z_1(s;t)+v_2 \tilde Z_2(s;t)   )  +v_2\tilde c^*_2(s )   \|^2 ds \nonumber \\
&= \left(\inf_{s \in [t, T]}\tilde Z_0(s)^2\right)\int_{t}^T   \| c_3(s) H^v(s) +v_1  c_3(s ) +v_2\tilde c^*_2(s )   \|^2 ds,
\end{align}
where $H^v(s):=v_1 (\tilde Z_1(s;t)-1 )+v_2 \tilde Z_2(s;t)  $ satisfies
\begin{align}\label{Hv dynamic}
dH^v(s)=\left( H^v(s)c_3(s)+v_1  c_3(s ) +v_2\tilde c^*_2(s )  \right)\tran dW(s), \quad s\in [t, T],~ H^v(t)=0.
\end{align}
Consequently, for any $v\in \bbR^2$ with $\|v\|=1$, $\eta>0$, and $p\ge 2$,
\begin{align}\label{Malliavin matrix probability}
&\prob( v\tran \Gamma v \leq \eta) \leq \prob\left( \left(\inf_{s \in [t, T]}\tilde Z_0(s)^2\right)\int_{t}^T   \| c_3(s) H^v(s) +v_1  c_3(s ) +v_2\tilde c^*_2(s )   \|^2 ds \leq  \eta \right) \nonumber\\
=& \prob\left(    \frac{\sup_{s \in [t, T]}\tilde Z_0(s)^{-2}}{ \int_{t}^T   \| c_3(s) H^v(s) +v_1  c_3(s ) +v_2\tilde c^*_2(s )   \|^2 ds  } \ge  1/\eta\right)\nonumber\\
 \le& \expect\left[\left|\frac{\sup_{s \in [t, T]}\tilde Z_0(s)^{-2}}{ \int_{t}^T   \| c_3(s) H^v(s)+v_1  c_3(s )  +v_2\tilde c^*_2(s )   \|^2 ds  } \right|^q\right]\big(1/\eta\big)^{-q}\nonumber\\
 \leq& \expect\left[\sup_{s \in [t, T]}\tilde Z_0(s)^{-4q}\right]^{1/2} \expect\left[ \left|\frac{1}{   \int_{t}^T   \| c_3(s) H^v(s)+v_1  c_3(s )  +v_2\tilde c^*_2(s )   \|^2 ds       }\right| ^{2q} \right]^{1/2}  \eta^{ q},
\end{align}
where we set $1/0=+\infty$.

Note that $\tilde Z_0(s) = \tilde Z_0(t)/ M(s)$, where
\begin{align*}
dM(s) = M(s) c_3(s)\tran  dW(s),\; s\in[t,T],\quad M(t)=1.
\end{align*}
As a result,
\begin{align}\label{first part}
& \expect\left[\sup_{s \in [t, T]}\tilde Z_0(s)^{-4q}\right]= \expect\left[\tilde Z_0(t)^{-4q}\left(\sup_{s \in [t, T]}M(s)\right)^{4q}\right]  \nonumber\\
&\le  \expect\left[\tilde Z_0(t)^{-8q}\right]^{1/2}\expect\left[\left(\sup_{s \in [t, T]}M(s)\right)^{8q}\right]^{1/2}<+\infty,
\end{align}
where the  last inequality is because  $c_3$ is piece-wise continuous and thus  bounded.
Thus, recalling \eqref{Malliavin matrix probability},   to prove \eqref{Malliavin matrix Lemma}, we only need to show for any $m\geq 2$,
\begin{align}\label{second part}
& \sup_{\|v\|=1} \expect\left[ | \frac{1}{   \int_{t}^T   \| c_3(s) H^v(s)+v_1  c_3(s )  +v_2\tilde c^*_2(s )   \|^2 ds      }|^{m}  \right]<+\infty.
\end{align}

For each $v=(v_1,v_2)\tran\in \bbR^2$ with $\|v\|=1$, because there exists $s\in [t, T)$, such that $v_1 c_3(s)+v_2 \tilde c^*_2(s) \neq 0$ and because $\tilde c^*_2$ and $c_3$ are right-continuous, we conclude that $\varphi(v):=\int_{t}^T \| v_1 c_3(s)+v_2\tilde c^*_2(s)  \|^2  d s >0$. It is straightforward to see that $\varphi(v)$ is continuous in $v$, so $L:=\inf_{\|v\|=1}\varphi(v)>0$. Because $\tilde c^*_2$ and $c_3$ are in $\mathfrak{C}_{\mathrm{pw}}([0, T);\mathbb{R}^d)$, they are bounded on $[0,T)$ by certain constant $\bar C>0$. Then, for any $v=(v_1,v_2)\tran\in \bbR^2$ with $\|v\|=1$,
\begin{align}\label{lower bound for quadratic variation}
& \int_t^T \| c_3(s) H^v(s)+v_1  c_3(s )  +v_2\tilde c^*_2(s )   \|^2 ds \nonumber\\
&=\int_t^T \| v_1  c_3(s )  +v_2\tilde c^*_2(s )   \|^2 ds  + \int_t^T \| c_3(s) H^v(s)\|^2 ds \nonumber\\
&\quad +2 \int_t^T H^v(s) c_3(s)\tran (v_1  c_3(s )  +v_2\tilde c^*_2(s ))      ds \nonumber\\
&\geq  L -T\sup_{s\in [t, T]} |H^v(s)|^2 \bar C^2-4T\sup_{s\in [t, T]} |H^v(s)| \bar C^2.
\end{align}
Because  $L -T \delta^2 \bar C^2-4T \delta \bar C^2 \geq L/2$, $\forall \delta\in [0, \bar \delta]$ for certain $\bar \delta>0$, then we conclude from \eqref{lower bound for quadratic variation} that there exists $\bar\epsilon>0$ such that for any $\epsilon\in (0, \bar\epsilon)$ and $v=(v_1,v_2)\tran\in \bbR^2$ with $\|v\|=1$,
\begin{align}\label{Exponential inequality  part 1}
& \prob(\int_t^T \| c_3(s) H^v(s)+v_1  c_3(s )  +v_2\tilde c^*_2(s )   \|^2 ds <\epsilon, \sup_{s\in [t, T]} |H^v(s)|< \epsilon^{1/4}    ) =0.
\end{align}
Thus, we have for any $\epsilon\in (0, \bar\epsilon)$,   and any $v=(v_1,v_2)\in \bbR^2$ with $\|v\|=1$,
	\begin{align}\label{Exponential inequality}
	& \prob\left(\int_t^T \| c_3(s) H^v(s)+v_1  c_3(s )  +v_2\tilde c^*_2(s )   \|^2 ds <\epsilon   \right)   \nonumber \\
	&= \prob\left(\int_t^T \| c_3(s) H^v(s)+v_1  c_3(s )  +v_2\tilde c^*_2(s )   \|^2 ds <\epsilon, \sup_{s\in [t, T]} |H^v(s)| \geq \epsilon^{1/4}   \right)   \nonumber \\
&\leq  2 e^{-\frac{1}{2}\epsilon^{-1/2}   },
	\end{align}
where the equality is the case due to \eqref{Exponential inequality  part 1} and the inequality is the case due to \eqref{Hv dynamic} and to the inequality (A.5) in \citet{NualartD:06mc}. Sending $\epsilon$ to 0 in the above, we immediately derive that $\int_t^T \| c_3(s) H^v(s)+v_1  c_3(s )  +v_2\tilde c^*_2(s )   \|^2 ds>0$ almost surely. For any $m\ge 2$, denote 
\begin{align*}
 Q(y):=\prob\left( \left(\int_t^T \| c_3(s) H^v(s)+v_1  c_3(s )  +v_2\tilde c^*_2(s )   \|^2 ds\right)^{-m}>y \right),\quad y\ge 0.
\end{align*}
Then, we derive from \eqref{Exponential inequality} that
	\begin{align}\label{Exponential inequality CDF}
&  Q(y)\leq 2 e^{-\frac{1}{2}y^{1/(2m)}   }, \quad y\geq \bar \epsilon^{-m}.
\end{align}
As a result,
\begin{align*}
  &\expect\left[\left(\int_t^T \| c_3(s) H^v(s)+v_1  c_3(s )  +v_2\tilde c^*_2(s )   \|^2 ds\right)^{-m}\right]\le \bar \epsilon^{-m} + \\
  &\quad +\expect\left[\left(\int_t^T \| c_3(s) H^v(s)+v_1  c_3(s )  +v_2\tilde c^*_2(s )   \|^2 ds\right)^{-m}\mathbf{1}_{\int_t^T \| c_3(s) H^v(s)+v_1  c_3(s )  +v_2\tilde c^*_2(s )   \|^2 ds< \bar \epsilon}\right]\\
&=\bar \epsilon^{-m}+ \int_{\bar \epsilon^{-m}}^\infty x d (1-Q(x))=\bar \epsilon^{-m}+\bar \epsilon^{-m} Q(\bar \epsilon^{-m} ) +\int_{\bar \epsilon^{-m}}^{+\infty }   Q(x) d x\\
  & \le \bar \epsilon^{-m}+2 e^{-\frac{1}{2}\bar \epsilon^{-1/2}   } \bar \epsilon^{-m} +\int_{\bar \epsilon^{-m}}^\infty 2 e^{-\frac{1}{2}x^{1/(2m)}} dx,
\end{align*}
where the second inequality is the case due to \eqref{Exponential inequality CDF}. Note that $\int_{\bar \epsilon^{-m}}^\infty 2 e^{-\frac{1}{2}x^{1/(2m)}} dx<+\infty$, so we immediately derive \eqref{second part}. The proof then completes.
  \halmos
\end{pfof}

\begin{pfof}{Lemma \ref{support set multiD}}
  \begin{enumerate}
  \item[(i)]
  For simplicity, we assume $\beta_0$, $\beta_{1}$, $\beta_{2,i}$, and $\beta_{3,i}$, $i=1,\dots, d$ are continuous on $[0,t_1)$ with the left-limit at $t_1$ existing and continuous on $[t_1,T]$ for certain $t_1\in (0,T)$. The case in which $[0,T]$ is divided into multiple pieces and the above functions are continuous on each of the pieces can be treated similarly.

  We first prove that ${\cal S}_{X^T}\subseteq \mathrm{cls}({\cal U})$. Consider any mollifier $\eta$, which is a non-negative, infinitely differentiable real function on $\bbR$ that is supported on $[0,1]$ and satisfies $\int_\bbR \eta(s)ds= 1$. For each $k\in \mathbb{N}$, define $W_k(t)=\big(W_{1,k}(t),...,W_{d,k}(t)\big)\tran$, where $W_{i,k}(t):=k\int_{\bbR} W_i(t-s)\eta(k s) ds$, $t\in [0, T]$ with the convention that $W_i(t)=0,t<0$. Then, $W_k\in \mathfrak{C}^\infty([0,T],\bbR^d)\subset \mathfrak{H}([0,T];\bbR^d)$. Because $\beta_0$, $\beta_{1}$, $\beta_{2,i}$, and $\beta_{3,i}$, $i=1,\dots, d$ are continuous on $[0,t_1)$ with the left-limit at $t_1$ existing, Theorem 3.1 of \citet{gyongy1989stability} yields that $\max_{t\in[0,t_1]}\|X(t)-f_{W_k}(t)\|$ converges to 0 in probability. In particular, $\|X(t_1)-f_{W_k}(t_1)\|$ converges to 0 in probability.

  Recall that $X(t) = \tilde X(t;X(t_1))$, $t\in [t_1,T]$, where $\tilde X(t,\xi)$ stands for the solution to the SDE \eqref{multi linear SDE} that starts at time $t_1$ with initial value $\xi$. We also have $f_{W_k}(t) = \tilde f_{W_k,f_{W_k}(t_1)}$, where $\tilde f_{w,x}$ stands for the solution to \eqref{x w multiD} that starts from time $t_1$ with initial value $x$. Consider $\tilde X(t;X(t_1)),t\in[t_1,T]$ on the filtered probability space $(\Omega, \cF, (\cF_t)_{t\in[t_1,T]},\prob)$ and recall $\|X(t_1)-f_{W_k}(t_1)\|$ converges to 0 in probability. Because $\beta_0$, $\beta_{1}$, $\beta_{2,i}$, and $\beta_{3,i}$, $i=1,\dots, d$ are continuous on $[t_1,T]$, Theorem 3.1 of \citet{gyongy1989stability} yields that $\max_{t\in[t_1,T]}\|\tilde X(t;X(t_1))-\tilde f_{W_k,f_{W_k}(t_1)}(t)\|$ converges in probability to 0. As a result, $\max_{t\in[t_1,T]}\|X(t)-f_{W_k}(t)\|$ converges in probability to 0 and, consequently, $\max_{t\in[0,T]}\|X(t)-f_{W_k}(t)\|$ converges in probability to 0. In other words, $f_{W_k}$, viewed as a random variable taking values on $\mathfrak{C}([0,T];\mathbb{R}^d$), converges in probability to $X$, viewed as a random variable on the same space. Then,
  \begin{align*}
    \prob(X\in \mathrm{cls}({\cal U})) \ge \limsup_{k\uparrow +\infty}\prob(f_{W_k}\in \mathrm{cls}({\cal U}))=1,
  \end{align*}
  showing that ${\cal S}_{X^T}\subseteq \mathrm{cls}({\cal U})$.

  Next, we prove ${\cal S}_{X^T}\supseteq \mathrm{cls}({\cal U})$. To this end, we first show that $\mathrm{cls}({\cal U}) = \mathrm{cls}(\bar {\cal U})$, where
  $\bar {\cal U}:= \{f_{w} \in \mathfrak{C}([0, T]; \bbR^l):   w\in \mathfrak{C}^\infty([0,T];\bbR^d) \}$. Because $\beta_0$, $\beta_{1}$, $\beta_{2,i}$, and $\beta_{3,i}$, $i=1,\dots, d$ are continuous on $[0,t_1)$ with the left-limit at $t_1$ existing, the first part of the proof of Theorem 4.1 in \citet{gyongy1989stability} yields that $\{f_{w} \in \mathfrak{C}([0, t_1]; \bbR^l):   w\in \mathfrak{C}^\infty([0,t_1];\bbR^d) \}$ is dense in $\{f_{w} \in \mathfrak{C}([0, t_1]; \bbR^l):   w\in \mathfrak{H}([0,t_1];\bbR^d) \}$. Similarly, $\{f_{w} \in \mathfrak{C}([t_1, T]; \bbR^l):   w\in \mathfrak{C}^\infty([t_1,T];\bbR^d) \}$ is dense in $\{f_{w} \in \mathfrak{C}([t_1,T]; \bbR^l):   w\in \mathfrak{H}([t_1,T];\bbR^d) \}$. Note that $f_w$ depends on $w$ via $w'$ only. As a result, $\mathrm{cls}({\cal U}) = \mathrm{cls}(\bar {\cal U})$. Therefore, in the following we only need to prove ${\cal S}_{X^T}\supseteq \mathrm{cls}(\bar {\cal U})$.

  Fix any $w \in \mathfrak{C}^\infty([0,T];\bbR^d)$ and $\epsilon>0$. Define $\tilde f_{w,x}(t),t\in[t_1,T]$ to be the solution to \eqref{x w multiD} that starts from time $t_1$ with initial value $x$. Then, it is straightforward to show that there exists $L>0$ such that
  \begin{align}\label{eq:LemmaSupportProof1}
  \max_{t\in[t_1,T]}\|\tilde f_{w,x}(t)-\tilde f_{w,y}(t)\|\le L\|x-y\|,\quad x,y\in \bbR^l.
  \end{align}
   Because $\beta_0$, $\beta_{1}$, $\beta_{2,i}$, and $\beta_{3,i}$, $i=1,\dots, d$ are continuous on $[0,t_1)$ with the left-limit at $t_1$ existing,
  Theorem 4.1 of \citet{gyongy1989stability} shows that ${\cal S}_{X^{t_1}}$ is the closure of  $\{f_{w} \in \mathfrak{C}([0, t_1]; \bbR^l):   w\in \mathfrak{H}([0,t_1];\bbR^d) \}$, so we have
  \begin{align}\label{eq:LemmaSupportProof2}
    \prob\left(\max_{t\in[0,t_1]}\|X(t) - f_{w}(t)\|<\frac{\epsilon}{3(L+1)}\right)>0.
  \end{align}
  Denote by $(\prob_{t_1,x})_{x\in \bbR^l}$ to be the family of probability measures conditional on $X(t_1) = x,x\in \bbR^l$. Then, because $\beta_0$, $\beta_{1}$, $\beta_{2,i}$, and $\beta_{3,i}$, $i=1,\dots, d$ are continuous on $[t_1,T]$, for each $x\in \bbR^l$, applying Theorem 4.1 of \citet{gyongy1989stability}, we conclude that the support of $\tilde X(t;x),t\in[t_1,T]$, which is the solution to the SDE \eqref{multi linear SDE} that starts at time $t_1$ with initial value $x$, is the same as the closure of $\{\tilde f_{w,x}  \in \mathfrak{C}([ t_1,T]; \bbR^l):   w\in \mathfrak{H}([t_1,T];\bbR^d) \}$, where $\tilde f_{w,x}$ is the solution to the ODE \eqref{x w multiD} that starts from time $t_1$ with initial value $x$. As a result,
  \begin{align*}
    \eta(x):=\prob_{t_1,x}\left( \max_{t\in[t_1,T]}\|\tilde X(t;x) - \tilde f_{w,x}(t)\|<\frac{\epsilon}{3} \right)>0.
  \end{align*}
  Then, denoting by $A$ the event that $\max_{t\in[0,t_1]}\|X(t) - f_{w}(t)\|<\frac{\epsilon}{3}$ and $\max_{t\in[t_1,T]}\|\tilde f_{w,f_w(t_1)}(t) -  \tilde f_{w,X(t_1)}(t)\|<\frac{\epsilon}{3}$, we have
  \begin{align*}
    &\prob\left(\max_{t\in[0,T]}\| X(t) -  f_{w}(t)\|<\epsilon \right)\ge  \prob\left(\max_{t\in[0,t_1]}\|X(t) - f_{w}(t)\|<\frac{\epsilon}{3}, \max_{t\in[t_1,T]}\|X(t) -  f_{w}(t)\|<\frac{2\epsilon}{3} \right)\\
    &= \prob\left(\max_{t\in[0,t_1]}\|X(t) - f_{w}(t)\|<\frac{\epsilon}{3}, \max_{t\in[t_1,T]}\|\tilde X(t;X(t_1)) -  \tilde f_{w,f_w(t_1)}(t)\|<\frac{2\epsilon}{3} \right)\\
    &\ge \prob\Big(A\text{ and } \max_{t\in[t_1,T]}\|\tilde X(t;X(t_1)) -  \tilde f_{w,X(t_1)}(t)\|<\frac{\epsilon}{3} \Big)\\
    & = \expect\Big[\mathbf 1_{A} \prob\left(\max_{t\in[t_1,T]}\|\tilde X(t;X(t_1)) -  \tilde f_{w,X(t_1)}(t)\|<\frac{\epsilon}{3}|\cF_{t_1}\right)\Big]\\
    & = \expect[\mathbf 1_{A}\eta(X(t_1))],
  \end{align*}
  where the last equality is the case due to the Markovian property of $X$. By \eqref{eq:LemmaSupportProof1} and \eqref{eq:LemmaSupportProof2}, we conclude that $\prob(A)>0$. Because $\eta(x)>0$ for all $x$, we conclude that
  \begin{align*}
    \prob\left(\max_{t\in[0,T]}\| X(t) -  f_{w}(t)\|<\epsilon \right)\ge  \expect[\mathbf 1_{A}\eta(X(t_1))]>0.
  \end{align*}
  This shows that ${\cal S}_{X^T}\supseteq \mathrm{cls}(\bar {\cal U})$.
  \item[(ii)] We first prove that $\mathrm{cls}(\mathbb{U}_t)\subseteq \mathbb{S}_{X(t)}$. For any $y\in \mathrm{cls}(\mathbb{U}_t)$, there exists $f_{n}\in {\cal U}$ such that $f_{n}(t)$ converges to $y$. Then, for each fixed $\epsilon>0$, there exists $n_\epsilon$ such that $\|f_{n_\epsilon}(t)-y\|<\epsilon/2$. Consequently,
      \begin{align*}
        \prob(\|X(t)-y\|<\epsilon)\ge \prob(\|X(t)-f_{n_\epsilon}(t)\|<\epsilon/2)\ge \prob(\max_{t\in[0,T]}\|X(t)-f_{n_\epsilon}(t)\|<\epsilon/2)>0,
      \end{align*}
      where the last inequality is the case because $f_{n_\epsilon}\in {\cal U}$ and ${\cal S}_{X^T}=\mathrm{cls}({\cal U})$.

      Next, we prove $\mathbb{S}_{X(t)}\subseteq \mathrm{cls}(\bar {\mathbb{U}}_t)$. To this end, consider any $y\notin\mathrm{cls}(\bar {\mathbb{U}}_t)$. Then, there exists $\epsilon_0>0$ such that $\|f(t)-y\|\ge \epsilon_0$ for any $f \in \bar {\cal U}$. Because ${\cal S}_{X^T}= \mathrm{cls}(\bar {\cal U})$, we have $\|f(t)-y\|\ge \epsilon_0$ for any $f \in {\cal S}_{X^T}$. Consequently,    $A:=\{f\in\mathfrak{C}([0, T]; \bbR^l)   :\|f(t)-y\|<\epsilon_0\}\subseteq {\cal S}_{X^T}^c$,
      where ${\cal S}_{X^T}^c$ stands for the complement of ${\cal S}_{X^T}$. As a result,
      \begin{align*}
        \prob(\|X(t)-y\|<\epsilon_0) = \prob( X\in A)\le \prob(X\in{\cal S}_{X^T}^c)=0,
      \end{align*}
      where the last equality is the case because ${\cal S}_{X^T}$ is the support of $X$. Thus, $y\notin \mathbb{S}_{X(t)}$. In other words, $\mathbb{S}_{X(t)}\subseteq \mathrm{cls}(\bar {\mathbb{U}}_t)$. Recall that $\bar {\mathbb{U}}_t \subset  {\mathbb{U}}_t$ by their definition. Thus, we have $\mathbb{S}_{X(t)}=\mathrm{cls}({\mathbb{U}}_t)=\mathrm{cls}(\bar {\mathbb{U}}_t)$.
\item[(iii)] Fix any $u,v\in \mathfrak{C}^\infty([0, T];\bbR^d) $ and for any $\lambda\in[0,1]$, define $w_\lambda:=\lambda u + (1-\lambda) v$. By Gronwall's inequality, we conclude from equation \eqref{x w multiD} that $L_1:=\sup_{\lambda\in[0,1]}\max_{t\in[0,T]}f_{w_\lambda}(t)<+\infty$. Moreover, there exists $L_2>0$ such that $|w_{\lambda_1,i}-w_{\lambda_2,i}|\le L_2|\lambda_1-\lambda_2|,\forall\lambda_1,\lambda_2\in[0,1]$. Applying Gronwall's inequality again, we conclude from equation \eqref{x w multiD} that $f_{w_\lambda}$ is continuous in $\lambda$. Therefore,  $\bar {\cal U}$ and $\bar {\mathbb{U}}_t$, $t\in[0,T]$, are connected.

    When $\ell=1$, because $\bar {\mathbb{U}}_t$ is connected, it is an interval. As a result, $\mathbb{S}_{X(t)} = \mathrm{cls}(\bar {\mathbb{U}}_t)$ is also an interval.  \halmos
\end{enumerate}
\end{pfof}

\section{Proofs}
\subsection{Proof of Theorem \ref{th:Density}}
By \eqref{CDFandQuantileEquation}, parts (i)--(iii) are consequences of straightforward calculation, so we only need to prove (iv) in the following.

Recall $\tilde Z_1(T;t)$ and $\tilde Z_2(T;t)$ as defined in Lemma \ref{lemma:Density}. It is straightforward to see that
  \begin{align*}
    \tilde X^*(T;t,x) = x \tilde Z_1(T;t) + \tilde Z_2(T;t).
  \end{align*}
  Noting that $\tilde Z_1(T;t)>0$, we derive
  \begin{align*}
    F^*(t,x,y) = \int_0^{+\infty}\int_{-\infty}^{y-z_1x  }g(z_1, z_2; t)dz_2 dz_1, 
  \end{align*}
  where $g(z_1,z_2;t)$ stands for the density of $(\tilde Z_1(T;t),\tilde Z_2(T;t))$. By Lemma \ref{lemma:Density},
  \begin{align}\label{eq:JointDenGrowthCon2}
    \sup_{z_1,z_2}(z_1^2+z_2^2)^{k/2}\left|\frac{\partial^{i+j}g}{\partial z_1^i\partial z_2^j}(z_1,z_2;t)\right|<+\infty
  \end{align}
  for any $k,i,j\ge 0$, so the dominated convergence theorem yields that $F^*(t,x,y)$ is differentiable in $y$ and its derivative is
  \begin{align}\label{density of X}
    \frac{\partial F^*}{\partial y}(t,x,y) = \int_0^{+\infty}g(z_1, y-z_1x; t) dz_1,
  \end{align}
  and that $F^*(t,x,y)$ is differentiable in $x$ with derivative
  \begin{align}\label{F x}
    \frac{\partial F^*}{\partial x}(t,x,y) = -\int_0^{+\infty}z_1g(z_1, y-z_1x; t) dz_1.
  \end{align}
  Moreover, $ \frac{\partial F^*}{\partial y}(t,x,y)$ and $ \frac{\partial F^*}{\partial x}(t,x,y)$ are continuous in $(x,y)$, so $F^*(t,x,y)$ is differentiable in $(x,y)$. Similar arguments show that $F^*$ is infinitely differentiable in $(x,y)$ with
  \begin{align*}
    \frac{\partial^{i+j}F^*}{\partial x^i\partial y^j} (t,x,y) = \int_0^{+\infty} \frac{\partial^{i+j}}{\partial x^i\partial y^j} \left(\int_{-\infty}^{y-z_1x  }g(z_1, z_2; t)dz_2\right) dz_1.
  \end{align*}
  Moreover, $F^*(t,x,y)$ and its derivatives with respect to $x$ and $y$ of any order are bounded in $(x,y)\in \bbR^2$.

  Finally, recalling \eqref{CDFandQuantileEquation}, we complete the proof. \halmos

\subsection{Proof of Theorem \ref{th:GeneralDensity}}
We need to consider $F^*(t,x,y)$ in following discussion and then recall the transformation \eqref{CDFandQuantileEquation}.

  We prove (i) first, and we only need to consider the case in which $t_*<t^*$. Straightforward calculation yields that $F^*\in \mathfrak{C}^{1,\infty}\big([t_{i-1},t_{i})\times (\bbR^2\backslash\{(\xi,\xi)\})\big)$ and \eqref{eq:FKFormulaGM} holds. In particular,
  \begin{align*}
    \frac{\partial F^*}{\partial t}(t,x,y) =
    \begin{cases}
    \phi\left(\frac{\ln \left(\frac{y-\xi}{x-\xi}\right)-\bar a_t}{\bar b_{t}}\right)\mathbf 1_{y>\xi}\left[g_1(t) \ln \left(\frac{y-\xi}{x-\xi}\right) + g_2(t)\right], & y\in \bbR, x>\xi,\\
    0, & y\neq \xi, x=\xi,\\
    -\phi\left(\frac{\ln \left(\frac{\xi-y}{\xi-x}\right)-\bar a_t}{\bar b_{t}}\right)\mathbf 1_{y<\xi}\left[g_1(t) \ln \left(\frac{\xi-y}{\xi-x}\right) + g_2(t)\right], & y\in \bbR, x<\xi,
    \end{cases}
  \end{align*}
  where $\bar b_t, \bar a_t$ are given in Theorem \ref{th:Density}-(iii), and
  \begin{align*}
    g_1(t) = \frac{d}{dt}(1/\bar b_t),\quad  g_2(t) = \frac{d}{dt}(-\bar a_t/\bar b_t).
  \end{align*}
  By the definition of $t^*$ and $t_*$, we have $\bar b_t>0$ for any $t<t^*$, so for any $\tau\in [t_*,t^*)$, $g_1(t)$ and $g_2(t)$ are bounded in $t\in[t_*,\tau]$ and $\phi\big((z-\bar a_t)/\bar b_t\big)P(z)$ are bounded in $(t,z)\in [t_*,\tau]\times \bbR$ for any polynomial function $P(z)$. Thus,
  \begin{align*}
     \sup_{t\in [t_*,\tau],(x,y)\neq (\xi,\xi)}\left|F_t^*(t,x,y)\right|<+\infty.
  \end{align*}
  Similar calculation shows that for any $j,k\in \mathbb{N}\cup \{0\}$, we have
  \begin{align}
    &\frac{\partial^{1+j+k} F^*}{\partial t\partial x^j \partial y^k}(t,x,y) = \mathbf 1_{y>\xi}
    (x-\xi)^{-j}(y-\xi)^{-k}\notag\\
    &\times \sum_{n=0}^{j+k}h_n(t,\ln \left(\frac{y-\xi}{x-\xi}\right) )\phi^{(n)}\left(\frac{\ln \left(\frac{y-\xi}{x-\xi}\right)-\bar a_t}{\bar b_{t}}\right), & y\in \bbR, x>\xi,\label{eq:GMPartialD}
  \end{align}
where $\phi^{(n)}$ stands for the $n$-th derivative of $\phi$ and $h_n(t,z)$ is certain function of $(t,z)$ such that $\sup_{t\in[t_*,\tau]}|h_n(t,z)|< C_n (1+|z|)$, $z\in \bbR$ for certain positive constant $C_n$. $\frac{\partial^{1+j+k} F^*}{\partial t\partial x^j \partial y^k}(t,x,y)$ takes a similar form when $x<\xi$ and is 0 when $x=\xi $ and $y\neq \xi$. For each fixed $\delta>0$, $(x,y)\notin B_2(\xi,\delta)$ implies that
\begin{align*}
  |x-\xi|^{-1}\le \delta^{-1} \left(1+\left(\frac{y-\xi}{x-\xi}\right)^2\right)^{1/2},\quad |y-\xi|^{-1}\le \delta^{-1} \left(1+\left(\frac{y-\xi}{x-\xi}\right)^{-2}\right)^{1/2}.
\end{align*}
Together with \eqref{eq:GMPartialD} and noting that for any $m\in \bbR$, $(1+ |\ln z|)(1+z^m)|\phi^{(n)}((\ln z-\bar a_t)/\bar b_t)|$ is bounded in $(t,z)\in [t_*,\tau]\times (0,+\infty)$, we immediately conclude \eqref{eq:GMDerivativeBound} for $\ell=1$. The case $\ell=0$ can be treated similarly.

Next, we prove (ii), and we only need to consider the case in which $t_*>0$. By the definition of $t_*$ and $t^*$, we conclude that for each $\tau\in (t_{m-1},t_*)$, it is either the case in which $c_3(s)=0,\forall s\in [\tau,T)$ or the case in which for any $v=(v_1,v_2)\tran \in \bbR^2$ with $\|v\|=1$, there exists $s\in [\tau,T)$ such that $v_1 c_3(s)+v_2\tilde c_2^*(s)\neq 0$. Thus, Theorem \ref{th:Density} yields that $F^*(\tau,x,y)$ is infinitely differentiable in $(x,y)\in \bbR^2$ with bounded derivatives. Recall that
\begin{align*}
  F^*(t,x,y) = \expect[F^*(\tau,\tilde X^*(\tau;t,x),y)].
\end{align*}
Because $\tilde X(\tau;t,x)$ is infinitely differentiable in $x$ pathwisely, and the derivatives of any order have finite arbitrary order  moments, the dominated convergence theorem yields that for any $j,k\in \mathbb{N}_0$,
\begin{align*}
  \frac{\partial^{j+k}F^*}{\partial x^j\partial y^k}(t,x,y) = \expect\left[\frac{\partial^{j+k}}{\partial x^j\partial y^k}F^*(\tau,\tilde X^*(\tau;t,x),y)\right].
\end{align*}
Because the derivatives of $\tilde X^*(\tau;t,x)$ with respect to $x$ is continuous in $(t,x)$, the dominated convergence theorem yields that $ \frac{\partial^{j+k}F^*}{\partial x^j\partial y^k}(t,x,y)$ is continuous and bounded in $(t,x,y)\in [t_{m-1},\tau)\times \bbR^2$. Moreover, Chapter 5,   Theorem  6.1 \citet{Friedman2012:SDEandApplications} implies that for each fixed $y\in \bbR$, $F^*(t,x,y)$, as a function of $(t,x)$, belongs to $\mathfrak{C}^{1,2}([t_{m-1},\tau)\times \bbR)$ and
\begin{align}\label{Fyeman Kac}
  F_t^*(t,x,y) = - \frac{1}{2}\| \tilde c_2^*(t) + c_3(t)x \|^2F_{xx}^*(t,x,y), \quad (t,x,y)\in [t_{m-1},\tau)\times \bbR^2,
\end{align}
showing that $F_t^*$ exists and $F_t^*(t,x,y)$ is continuous in $(t,x,y)\in [t_{m-1},\tau)\times \bbR^2$ and infinitely differentiable in $(x,y)\in \bbR^2$, and $\sup_{t\in[t_{m-1},\tau),y\in \bbR}|\frac{\partial^{j+k}F_t^*}{\partial x^j\partial y^k}(t,x,y)|$ is of polynomial growth in $x$ for any $j,k\in \mathbb{N}_0$. Because $\tau$ is arbitrarily, we conclude that $F^*\in \mathfrak{C}^{1,\infty}\big([t_{m-1},t_*)\times \bbR^2\big)$, $F_t^*\in \mathfrak{C}^{0,\infty}\big([t_{m-1},t_*)\times \bbR^2\big)$, and for any $\tau'\in[t_{m-1},t_*)$ and $j,k\in \mathbb{N}_0$, $|\frac{\partial^{j+k}F^*}{\partial x^j\partial y^k}(t,x,y)|$ is bounded in $(t,x,y)\in [t_{m-1},\tau']\times \bbR^2$ and $\sup_{t\in[t_{m-1},\tau'],y\in \bbR}|\frac{\partial^{j+k}F_t^*}{\partial x^j\partial y^k}(t,x,y)|$ is of polynomial growth in $x$. Similar arguments show that for any $i=1,\dots, m-1$, $F^*\in \mathfrak{C}^{1,\infty}\big([t_{i-1},t_{i})\times \bbR^2\big)$, $F_t^*\in \mathfrak{C}^{0,\infty}\big([t_{i-1},t_{i})\times \bbR^2\big)$, $|\frac{\partial^{i+j}F^*}{\partial x^i\partial y^j}(t,x,y)|$ is bounded in $(t,x,y)\in [t_{i-1},t_{i})\times \bbR^2$, and $\sup_{t\in[t_{i-1},t_{i}),y\in \bbR}|\frac{\partial^{i+j}F_t^*}{\partial x^i\partial y^j}(t,x,y)|$ is of polynomial growth in $x$.
Combining \eqref{tilde c2 *}, \eqref{CDFandQuantileEquation}, and \eqref{Fyeman Kac}, we derive that ${\cal A}F(t,x,y)=0$.

Next, we prove (iii), and we only need to consider the case $t_*\in(0,t^*)$. For any $t\in [0,t_*)$, we have $F^*(t,x,y) = \expect[F^*(t_*,\tilde X^*(t_*;t,x),y)]$. Fixing $(x_0,y_0)\neq (\xi,\xi)$ and we prove that $F^*(t,x,y)$ is continuous at $(t_*,x_0,y_0)$ from the left of $t_*$.
When $y_0\neq \xi$, because $F^*(t_*,x',y')$ is continuous in $(x',y')\neq (\xi,\xi)$ and because $\tilde X^*(t_*;t,x)$ is continuous in $(t,x)$ pointwisely, the dominated convergence theorem yields that $\lim_{t\uparrow t_*,(x,y)\rightarrow (x_0,y_0)}F^*(t,x,y) = F^*(t_*,x_0,y_0)$. When $x_0\neq  \xi$ and $y_0=\xi$, Corollary \ref{coro:density} shows that $\prob(\tilde X^*(t_*;t,x)=\xi)=0$ for any $t$ and $x\neq \xi$, so the dominated convergence theorem again yields that $\lim_{t\uparrow t_*,(x,y)\rightarrow (x_0,y_0)}F^*(t,x,y) = F^*(t_*,x_0,y_0)$.

The same calculation as in the proof of part (ii) of the theorem yields that for any $j,k\in \mathbb{N}_0$
\begin{align}\label{eq:GMDerivativeBound2}
\sup_{x\in \bbR}\left|\frac{\partial^{j+k}F^*}{\partial x^j\partial y^k}(t_*,x,y)\right|\le C_{j,k}|y-\xi|^{-(j+k)},\quad y\neq \xi
 \end{align}
 for some constant $C_{j,k}>0$.
 Because $\tilde X(t_*;t,x)$ is infinitely differentiable in $x$ pathwisely, and the derivatives of any order have finite moments of any order, the dominated convergence theorem immediately yields that
\begin{align*}
  \frac{\partial^{j+k}F^*}{\partial x^j\partial y^k}(t,x,y) = \expect\left[\frac{\partial^{j+k}}{\partial x^j\partial y^k}F^*(t_*,\tilde X^*(t_*;t,x),y)\right],\quad t\in[0, t_*],x\in\bbR,y\neq \xi
\end{align*}
and that $\frac{\partial^{j+k}F^*}{\partial x^j\partial y^k}$ is continuous at $(t_*,x,y)$ from the left of $t_*$ for any $x\in \bbR$ and $y\neq \xi$.
 Moreover, \eqref{eq:GMDerivativeBound2} yields that for any $\delta>0$, $\sup_{t\in [0,t_*],x\in \bbR, |y-\xi|>\delta}\left|\frac{\partial^{j+k} F^*}{\partial x^j \partial y^k}(t,x,y)\right|<+\infty$, which together with \eqref{eq:GMDerivativeBound}  yields that \eqref{eq:GMDerivativeBoundGlobal} holds.
  Recalling  \eqref{eq:FKFormulaNonGM}, we have $\sup_{t\in [0,\tau],|y-\xi|>\delta}\left|\frac{\partial^{j+k} F_t^*}{\partial x^j \partial y^k}(t,x,y)\right|$ is of polynomial growth in $x$ for any $\delta>0$ and $\tau\in [0,t^*)$.

Now, we prove (iv), and only need to consider the case $t_*>0$.
 For each $t\in [0, t_*)$, it is either the case in which $c_3(s)=0,\forall s\in [t,T)$ or the case in which for any $v=(v_1,v_2)\tran \in \bbR^2$ with $\|v\|=1$, there exists $s\in [t,T)$ such that $v_1c_3(s)+v_2 \tilde c_2^*(s)\neq 0$.
 When $c_3(s)=0,\forall s\in [t,T)$ , the definition of $t_*$ implies that $c_2(s)\neq 0$ for some $s\in [t, t_*)$, and a direct calculation shows that $\frac{\partial F^*}{\partial x}(t,x,y)<0$, for any $(x,y)\in \bbR^2$.
 For the later case and for sake of contradiction, we assume there exists  $(x,y)\in  \bbR^2$, such that $F^*(t,x,y)\in (0, 1)$ and  $\frac{\partial F^*}{\partial x}(t,x,y)\geq 0$.
 According to \eqref{F x}, we obtain that $\frac{\partial F^*}{\partial x}(t,x,y)=0$, which together with \eqref{density of X} implies that $ \frac{\partial F^*}{\partial y}(t,x,y)=0$. Theorem \ref{th:Density} yields that $\tilde X^*(T;t,x)$ has a continuous density and Lemma \ref{support set multiD}-(iii) then implies that the interior of the support of $\tilde X^*(T;t,x)$ is a nonempty interval, so $y$ is in the interior of the support from $F^*(t,x,y)\in (0, 1)$. Proposition 2.1.8 of \citet{NualartD:06mc} further implies that the density, i.e. $\frac{\partial F^*}{\partial y}(t,x,y)$ is positive when $y$ is in the interior of the support of $\tilde X^*(T;t,x)$. Then, we arrive a contradiction.

 When $t_*<t^*$, for any $t\in [t_*, t^*)$, straightforward calculation yields that
  \begin{align*}
    \frac{\partial F^*}{\partial x}(t,x,y) =
    \begin{cases}
    -\phi\left(\frac{\ln \left(\frac{y-\xi}{x-\xi}\right)-\bar a_t}{\bar b_{t}}\right)(\mathbf 1_{x>\xi})\frac{1}{\bar b_{t}(x-\xi) }, & x\in \bbR, y>\xi\\
    0, & x\neq \xi, y=\xi,\\
     \phi\left(\frac{\ln \left(\frac{\xi-y}{\xi-x}\right)-\bar a_t}{\bar b_{t}}\right)(\mathbf 1_{x<\xi})\frac{1}{\bar b_{t}(x-\xi) }, & x\in \bbR, y<\xi,
    \end{cases}
  \end{align*}
Because for fixed $t\in [t_*, t^*)$, $F^*(t,x,y)\in (0, 1)$ if and only if $x>\xi, y>\xi$ or $x<\xi, y<\xi$, then $\frac{\partial F^*}{\partial x}(t,x,y)<0$.

Finally, we prove (v). We consider the case $x>y$ only, as the case $x<y$ can be treated similarly. Because $ \tilde c_2^*,c_3\in \mathfrak{C}_{\mathrm{pw}}([0, T];\mathbb{R}^d)$, there exists $C>0$ such that $\expect[|\tilde X^*(t^*;t,x')-x'|^{2}]\le C|t^*-t|$ for all $x'$ and $t$ that is sufficiently close to $t^*$. By Chebyshev's inequality, for all $t$ that is sufficiently close to $t^*$ and all $x'>y'$,
\begin{align*}
  F^*(t,x',y') &= \prob(\tilde X^*(t^*;t,x')\le y') \le  \prob(|\tilde X^*(t^*;t,x')-x'|\ge x'-y')\\
  &\le \expect[|\tilde X^*(t^*;t,x')-x'|^{2}]/(x'-y')^2\le C|t^*-t|/(x'-y')^2,
\end{align*}
which immediately implies $\lim_{t\uparrow t^*,(x',y')\rightarrow (x,y)}F^*(t,x',y')=0=F^*(t^*,x,y)$. \halmos

\subsection{Proof of Corollary \ref{coro:GeneralQuantiles}}

Because of \eqref{CDFandQuantileEquation}, we need to consider $G^*(t,x,y)$ in following discussion.
  (i) is trivial to prove, so next we assume $t^*>0$ and prove (ii) and (iii).

  For any $t\in[t_*,t^*)$ and $x\neq \xi$, $\tilde X^*(T;t,x)$ has a continuous density, the support of $\tilde X^*(T;t,x)$ is an interval, and the density is positive in the interior of the support of $\tilde X^*(T;t,x)$. Thus, $G^*(t,x,\alpha)$ is uniquely determined by $ F^*(t,x,G^*(t,x,\alpha)) = \alpha$, and $F^*_y(t,x,G^*(t,x,\alpha))>0$. Suppose $t_*>0$ and consider any $t\in [0,t_*)$ and $x\in \bbR$. Because $t<t_*$, Theorem \ref{th:Density} yields that $\tilde X^*(T;t,x)$ has a continuous density and Lemma \ref{support set multiD}-(iii) then implies that the interior of the support of $\tilde X^*(T;t,x)$ is a nonempty interval. Proposition 2.1.8 of \citet{NualartD:06mc} further implies that the density is positive in the interior of the support of $\tilde X^*(T;t,x)$. Thus, $G^*(t,x,\alpha)$ is uniquely determined by $ F^*(t,x,G^*(t,x,\alpha)) = \alpha$ and $F^*_y(t,x,G^*(t,x,\alpha))>0$. Then, (ii) and (iii) of the corollary just follows from the implicit function theorem and Theorem \ref{th:GeneralDensity}.

  Finally, we prove (iv). Fix any $x\in \bbR$ and $\alpha \in (0,1)$. For the sake of contradiction, suppose $G^*(t,x',\alpha')$ does not converge to $x$ when $(t,x',\alpha')\rightarrow (t^*,x,\alpha)$. Then, there exists a sequence $(t_n,x_n,\alpha_n)$ that converges to $(t^*,x,\alpha)$ and satisfies either $G^*(t_n,x_n,\alpha_n)\ge x+\delta,n\in\mathbb{N}$ or $G^*(t_n,x_n,\alpha_n)\le x-\delta,n\in\mathbb{N}$ for some $\delta>0$. If $t_*<t^*$ and $x_n=\xi$ for infinitely many $n$, we immediately have $x=\xi$ and $G^*(t_n,x_n,\alpha_n) = G^*(t_n,\xi,\alpha_n)=\xi$ for certain $n$, which is a contradiction. If $t_*=t^*$ or  $x_n=\xi$ only for finitely many $n$, we conclude that for sufficiently large $n$, $F^*(t_n,x_n,G^*(t_n,x_n,\alpha_n))=\alpha_n$. In the case in which $G^*(t_n,x_n,\alpha_n)\ge x+\delta,n\in\mathbb{N}$, we have $\alpha_n=F^*(t_n,x_n,G^*(t_n,x_n,\alpha_n))\ge F^*(t_n,x_n,x+\delta)$. Because $\alpha_n\rightarrow \alpha\in (0,1)$ and $F^*(t_n,x_n,x+\delta)\rightarrow F^*(t^*,x,x+\delta)=1$ by Theorem \ref{th:GeneralDensity}-(v), we arrive at contradiction. Similarly, in the case in which $G^*(t_n,x_n,\alpha_n)\le x-\delta,n\in\mathbb{N}$, we can also derive contradiction. The proof then completes.    \halmos

\subsection{Proof of Theorem \ref{th:GenralSupport}}

Because of the transformation \eqref{MappingSupport}, we only need to consider the support of the distribution of $X^*(t)$, denoted by $\mathbb{S}_{ X^*(t)}$, and the set of  reachable states of $X^*$ at time $t$, denoted by $\bbX^*_t$. For any $t\in [0, \underline{t})$, we have $c_2^*(t)=0$ and thus $h(t)=0$. As a result, $\underline{t} \leq \bar{t}$. For readability, we divide the remaining proof into several steps.

\subsubsection{Characterize $\mathbb{S}_{ X^*(t)}$ by an optimal control problem}\label{subsubse:PfSupportChar}


Lemma \ref{support set multiD} shows that $\mathbb{S}_{X^*(t)}$ is a closed interval with the lower end $\underline{x}^*(t)\in[-\infty,+\infty)$ and the upper end $\bar{x}^*(t)\in(-\infty,+\infty]$. Moreover,
  \begin{align*}
  &\underline{x}^*(t) = \inf_{w\in \mathfrak{H}([0,T];\bbR^d)}f_{w}(t) = \inf_{w\in \mathfrak{C}^\infty([0,T];\bbR^d)}f_{w}(t)=\inf_{w\in \hat{\mathfrak{H}}([0,T];\bbR^d)}f_{w}(t),\\
  &\bar{x}^*(t) = \sup_{w\in \mathfrak{H}([0,T];\bbR^d)}f_{w}(t) = \sup_{w\in \mathfrak{C}^\infty([0,T];\bbR^d)}f_{w}(t)=\sup_{w\in \hat{\mathfrak{H}}([0,T];\bbR^d)}f_{w}(t),
\end{align*}
where $\hat{\mathfrak{H}}([0,T];\bbR^d)$ denotes the set of absolutely continuous $w$ with a bounded derivative and for each $w\in \mathfrak{H}([0,T];\bbR^d)$,
$f_w$ is the solution to the following ODE:
  \begin{align}\label{f w one-dimensional}
f_{w}'(t)&= H^w(t,f_{w}(t)), t\in [0, T], \quad f_{w}(0)=0.
\end{align}
Here,
\begin{align}
  H^w(t,x):= \big(c^*_2(t)+c_3(t)x\big)\tran \big(w'(t)-\frac{1}{2}c_3(t)\big),\quad t\in [0, T],x\in \bbR.
\end{align}
Then, straightforward calculation yields
  \begin{align*}
f_w(t)&=\int_0^t e^{\int_s^t  c_3(\tau)\tran \left( -\frac{1}{2}c_3(\tau) +  w'(\tau)\right)d\tau  }c_2^*(s)\tran \big( -\frac{1}{2}c_3(s) +  w'(s)\big) ds,\quad t\in [0, T].
\end{align*}
As aresult, for any $t<s$, by considering a particular $w\in \mathfrak{H}([0,T];\bbR^d)$ with $w'(\tau)=\frac{1}{2}c_3(\tau), \tau\in (t, s] $, it is straightforward to see that
\begin{align*}
  \inf_{w\in \mathfrak{H}([0,T];\bbR^d)}f_{w}(s)\le \inf_{w\in \mathfrak{H}([0,T];\bbR^d)}f_{w}(t),\quad \sup_{w\in \mathfrak{H}([0,T];\bbR^d)}f_{w}(s)\ge \sup_{w\in \mathfrak{H}([0,T];\bbR^d)}f_{w}(t),
\end{align*}
so $\underline{x}^*(s)\le \underline{x}^*(t)$ and $\bar{x}^*(s)\ge \bar{x}^*(t)$, i.e., $\mathbb{S}_{X^*(t)}\subseteq \mathbb{S}_{X^*(s)}$.

Next, we prove that $\underline{x}^*$ and $\bar{x}^*$ are left-continuous. Fix any $t\in (0,T]$. Because $\underline{x}^*$ is decreasing, we conclude $\underline{x}^*(t)\le \liminf_{s\uparrow t}\underline{x}^*(s)$. On the other hand, for any $w\in  \hat{\mathfrak{H}}([0,T];\bbR^d)$, we have
\begin{align*}
  f_w(t) = \lim_{s\uparrow t}f_w(s) \ge \limsup_{s\uparrow t}\underline{x}^*(s).
\end{align*}
As a result, $\underline{x}^*(t) = \inf_{w\in \hat{\mathfrak{H}}([0,T];\bbR^d)}f_{w}(t) \ge \limsup_{s\uparrow t}\underline{x}^*(s)$. Thus, $\underline{x}^*(t) = \lim_{s\uparrow t}\underline{x}^*(s)$. Similarly, $\bar{x}^*(t) = \lim_{s\uparrow t}\bar{x}^*(s)$.

Because $\mathbb{S}_{X(t)} = \lambda_0(t) + \lambda_1(t) \left(x_0 + \mathbb{S}_{X^*(t)} \right)$, we derive $\underline{x}(t) = \lambda_0(t) + \lambda_1(t) \left(x_0 + \underline{x}^*(t) \right)$ and $\bar{x}(t) = \lambda_0(t) + \lambda_1(t) \left(x_0 + \bar{x}^*(t) \right)$. Because $\lambda_0$ and $\lambda_1$ are continuous, we conclude that $\underline{x}$ and $\bar{x}$ are left-continuous.

Next, we solve $\underline{x}^*(t)$ and $\bar{x}^*(t)$. For any $w\in \mathfrak{\hat H}([0,T];\bbR^d)$ and $(t,x)\in [0, T]\times \bbR$, define $V^{w}(t,x;\tau),\tau\in [t,T]$ to be the solution to the following equation
\begin{align*}
  \frac{\partial V^{w}}{\partial \tau}(t,x;\tau) = H^w\big(\tau, V^{w}(t,x;\tau)\big),\quad \tau \in [t,T],\quad V^{w}(t,x;t)=x.
\end{align*}
By definition, $f_w(s) = V^{w}(0,0;s),s\in[0,T]$. Straightforward calculation shows that
  \begin{align}\label{f w}
&V^{w}(t,x;\tau):=x e^{\int_{t}^{\tau} c_3(z)\tran \left( -\frac{1}{2}c_3(z) +  w'(z)\right) dz   } \nonumber \\
&\quad +\int_{t}^\tau e^{\int_s^\tau c_3(z)\tran \left( -\frac{1}{2}c_3(z) +  w'(z)\right) dz   }c_2^*(s)\tran \big( -\frac{1}{2}c_3(s) +  w'(s)\big)ds. 
\end{align}
Consequently, denoting by $V_t^{w}$ and $V_x^{w}$ the partial derivatives of $V^w$ with respect to $t$ and $x$, respectively, we derive, for any fixed $\tau\in (0,T]$, that
  \begin{align}\label{ODEFeymanKac}
V_t^{w}(t, x;\tau )+V_x^{w}(t,x;\tau )H^w(t,x) = 0, \quad \forall x\in \bbR\text{ and almost everywhere } t\in [0, \tau].
\end{align}
As a result, for any $w, \hat w \in \mathfrak{\hat H}([0,T];\bbR^d)$ and $\tau\in (0, T]$, we have
     \begin{align}\label{DifferenceOfTerminalValue}
&f_{\hat w}(\tau)-f_{w}(\tau)=V^{w}(\tau,f_{\hat w}(\tau);\tau)-V^{w}(0,0;\tau)=\int_{0}^\tau \frac{\partial V^{w}}{\partial s}(s, f_{\hat w}(s) ;\tau)  ds \nonumber \\
&=\int_{0}^\tau \big(V_t^{w}(s, f_{\hat w}(s);\tau )+V_x^{w}(s, f_{\hat w}(s);\tau )H^{\hat w}(s,f_{\hat w}(s)) \big)ds \nonumber \\
&=\int_{0}^\tau V_x^{w}(s, f_{\hat w}(s);\tau )\big(H^{\hat w}(s,f_{\hat w}(s)) - H^{ w}(s,f_{\hat w}(s)) \big)ds \nonumber \\
&=\int_{0}^\tau e^{\int_s^\tau c_3(z)\tran \left( -\frac{1}{2}c_3(z) +  w'(z)\right)dz   } \big( c_2^*(s)+c_3(s)f_{\hat w}(s)\big)\tran\big(\hat w'(s)-w'(s) \big)ds.
\end{align}

Because $c_3$ is in $\mathfrak{C}_{\mathrm{pw}}([0, T];\bbR^d)$ and thus bounded on $[0,T]$, there exists a constant $L>0$ such that
\begin{align*}
  \big|\|c_2^*(s)+c_3(s)x\|-\|c_2^*(s)+c_3(s)y\|\big|\le \|c_3(s)x-c_3(s)y\|\le L|x-y|,\quad \forall x,y\in \bbR,s\in [0,T].
\end{align*}
As a result, for any constant $k\in \bbR$, we can define $g_k$ to be the unique solution to the following equation:
 \begin{align}\label{IntegralEquation}
g_k(t)=\int_0^t k\|c_2^*(s)+c_3(s)g_k(s)\|ds,\quad t\in [0, T].
\end{align}
Define $w_k(t):=\int_0^t \frac{c_3(s)}{2}ds+k \int_{0}^t  \frac{ c_2^*(s)+c_3(s)g_k(s)  }{\| c_2^*(s)+c_3(s)g_k(s)  \| } \mathbf{1}_{E_k}(s)   ds$, $t\in [0, T]$, where $E_k$ is the set of $s\in[0,T]$ such that $c_2^*(s)+c_3(s)g_k(s)\neq 0$. Then, $w_k\in \mathfrak{\hat H}([0,T];\bbR^d)$. Straightforward calculation from \eqref{f w one-dimensional} yields
    \begin{align}\label{IntegralEquationCopy}
f_{w_k}(t)=\int_0^t k \frac{[c_2^*(s)+c_3(s)f_{w_k}(s) ]\tran[c_2^*(s)+c_3(s)g_k(s) ] }{\| c_2^*(s)+c_3(s)g_k(s)  \| } \mathbf{1}_{E_k}(s)ds,\quad t\in [0,T].
\end{align}
Comparing \eqref{IntegralEquation} and \eqref{IntegralEquationCopy}, we derive $f_{w_k}=g_k$. Moreover, by the standard comparison theorem, $g_k$ is increasing in $k$.

Now, taking $\hat w=w_k$ in \eqref{DifferenceOfTerminalValue} and recalling $g_k=f_{w_k}$, $t\in [0, T]$, we derive
  \begin{align}\label{DifferenceOfTerminalValue w_k}
&g_k(\tau)-f_{w}(\tau)=\int_{0}^\tau e^{\int_s^\tau c_3(z)\tran \big( -\frac{1}{2}c_3(z) +  w'(z)\big) dz   } [ c_2^*(s)+c_3(s)g_k(s)]\tran[w_k'(s)-w'(s) ]ds \nonumber\\
&=k\int_{0}^\tau e^{\int_s^\tau c_3(z)\tran \left( -\frac{1}{2}c_3(z) +  w'(z)\right) dz   }\|c_2^*(s)+c_3(s)g_k(s)  \| ds \nonumber \\
&\quad -\int_{0}^\tau e^{\int_s^\tau c_3(z)\tran \left( -\frac{1}{2}c_3(z) +  w'(z)\right) dz   }\big(c_2^*(s)+c_3(s)g_k(s)\big)\tran\big(w'(s)-\frac{1}{2}c_3(s)\big)ds.
\end{align}
Denote $G^{w}_k(\tau):=\int_{0}^\tau e^{\int_s^\tau c_3(z)\tran\left( -\frac{1}{2}c_3(z) +  w'(z)\right) dz   }\|c_2^*(s)+c_3(s)g_k(s)  \| ds$.
 Because $c_3(t)$ and $w'(t)$ are bounded for a.e. $t\in [0, T]$, we conclude from \eqref{DifferenceOfTerminalValue w_k} that there exists constant $L_w>0$, such that
  \begin{align}
(k-L_w)G^{w}_k(\tau)\le g_k(\tau)-f_{w}(\tau)\le (k+L_w)G^{w}_k(\tau),\quad \forall k\in \bbR. \label{DifferenceOfTerminalValue w_k 1}
\end{align}
Because \eqref{DifferenceOfTerminalValue w_k 1} is true for all $w \in \mathfrak{\hat H}([0,T];\bbR^d)$, we immediately derive that
 \begin{align}\label{UpperAndLowerBoundOfSupport}
\bar x^*(t) =\sup_{w\in \hat{\mathfrak{H}}([0,T];\bbR^d)}f_{w}(t) =  \lim_{k\uparrow +\infty}g_k(t), \quad \underline{x}^*(t) =\inf_{w\in \hat{\mathfrak{H}}([0,T];\bbR^d)}f_{w}(t) =  \lim_{k\downarrow -\infty}g_k(t).
\end{align}

\subsubsection{Proof of Part (i)}
By the definition of $\underline{t}$, we conclude $X^*(t)=0$, $t\in[0,\underline{t}]$ and thus $X(t)=\lambda_0(t)+\lambda_1(t) x_0$, $t\in[0,\underline{t}]$, so  part (i) of the theorem follows immediately.

\subsubsection{Proof of Part (ii)}\label{subsubse:PfSupportPartii}

Fix $t\in (\underline{t},T]$. Note that $X^*(t) = X^*(t;\underline{t},0)$. Recalling the definition of $\underline{t}$ and applying Corollary \ref{coro:density}-(ii), with $[t,T]$, $x$, and $\tilde c_2^*$ therein set to be $[\underline{t},t]$, 0, and $c_2^*$, respectively, we immediately conclude that
%
$X^*(t) = X^*(t;\underline{t},0)$ possesses a probability density function. Then, $X(t)=\lambda_0(t)+\lambda_1(t) (x_0+  X^*(t)) $ possesses a probability density function, $\mathrm{int}(\mathbb{S}_{X(t)})$ is a nonempty open interval and $\bbX_t=\mathrm{int}(\mathbb{S}_{X(t)})$.


Now, fix $t\in (\bar t,T]$. 
For each $s\in[0,T]$, it follows from the definition of $h$ that $\min_{ x\in \bbR} \| c_2^*(s)+ c_3(s)x\|=\|c_2^*(s)+c_3(s)h(s)\|$. As a result,
\begin{align*}
|g_k(t)| = |k|\int_0^t\|c_2^*(s)+c_3(s)g_k(s)\|ds \ge |k|\int_0^t\|c_2^*(s)+c_3(s)h(s)\|ds,
\end{align*}
where the first equality comes from \eqref{IntegralEquation}. Sending $k$ in the above to $+\infty$ and $-\infty$, respectively, and recalling that $\int_0^t\|c_2^*(z)+c_3(z)h(z)\|dz>0 $ because $t\in(\bar{t}, T]$, we immediately conclude that $\bar x^*(t) =\lim_{k\uparrow +\infty}g_k(t) =+\infty$ and $\underline{x}^*(t)=\lim_{k\downarrow -\infty}g_k(t)=-\infty$, i.e., $\bbX_t = \bbX_t^*=\bbR$.

\subsubsection{Proof of Part (iii)}
Recall that $g_k$ as defined in \eqref{IntegralEquation}. As already showed in Section \ref{subsubse:PfSupportChar}, $g_k$ is increasing in $k$. Moreover, $g_0(s)=0,\forall s\in [0,T]$, so $g_k(s)\ge 0,s\in[0,T]$ for all $k\ge 0$ and $g_k(s)\le 0,s\in [0,T]$ for all $k\le 0$. In addition, it is straightforward to see from \eqref{IntegralEquation} that for $k\ge 0$, $g_k(s)$ is increasing in $s\in [0,T]$ and for $k\le 0$, $g_k(s)$ is decreasing in $s\in [0,T]$.

Assume $\underline{t}<\bar{t}$ and fix $t\in (\underline{t},\bar{t}]$ in the following. We claim that $\tau_t>\underline{t}$. Otherwise, we have $c_3(z)=0,z\in (\underline{t},\bar t)$. By the definition of $\bar{t}$, we have $c_2(z)+c_3(z)h(z)=0$ for almost everywhere $z\in [\underline{t},\bar{t}]$. By the right-continuity of $c_2$, we derive $c_2(z)=0,\forall z\in [\underline{t},\bar t)$. Then, by the definition of $\underline{t}$, we have $\underline{t}\ge \bar {t}$, which is a contradiction. Thus, we must have $\tau_t>\underline{t}$.

By the definition of $\bar{t}$, we have $c_2^*(z)  = - c_3(z)h(z)$ for almost everywhere $z\in [0,\bar t)$. As a result, we conclude from \eqref{IntegralEquation} that
\begin{align}
  g_k(s) = k\left(\int_0^s\|c_3(z)\|| g_k(z)-h(z)|dz\right),\; s\in [0,\bar t].
\end{align}
By the definition of $\underline{t}$, $c_2^*(z)=0$ and thus $h(z)=0$ for all $z\in [0,\underline{t})$. As a result, $g_k(s) = k\left(\int_0^s\|c_3(z)\|| g_k(z)|dz\right),\forall s\in [0,\underline{t}]$, which implies $g_k(s) = 0,\forall s\in [0,\underline{t}]$, and
\begin{align}\label{IntegralEquation2}
  g_k(s) = k \int_{\underline{t}}^s\|c_3(z)\|| g_k(z)-h(z)|dz,\; s\in [\underline{t},\bar t].
\end{align}

Suppose that there exists $s\in [\underline{t},t)$ such that $h(s)<0$. Then, by the definition of $h$, we must have $c_3(s)\neq 0$. By the right-continuity of $c_2^*$ and $c_3$, there exists $\epsilon_0>0$ and $\delta\in (0,t-s)$ such that $h(z)\le -\epsilon_0$ and $\|c_3(z)\|\ge \epsilon_0$ for all $z\in [s,s+\delta]$. As a result, for any $k\ge 0$,
\begin{align*}
  \int_{\underline{t}}^t\|c_3(z)\| | g_k(z)-h(z)|dz\ge  \int_{s}^{s+\delta}\|c_3(z)\| | g_k(z)-h(z)|dz\ge \epsilon_0^2\delta,
\end{align*}
where the last inequality is the case because $g_k(z)\ge 0,z\in [0,T]$ for all $k\ge 0$. We then conclude from \eqref{IntegralEquation2} that $\bar{x}^*(t) = \lim_{k\uparrow +\infty}g_k(t) = +\infty$.

A similar argument shows that if there exists $s\in [\underline{t},t)$ such that $h(s)>0$, then $\underline{x}^*(t) = -\infty$. As a result, $\bbX_t = \bbX_t^*=\bbR$.

Next, we consider the case in which $h(s)\le 0$ for all $s\in [\underline{t}, t)$ and there exists $s_1,s_2\in [\underline{t},t)\cap D$ with $s_1<s_2$ such that $h(s_1)<h(s_2)$. Then, because $h(s_1)<0$, as shown in the above, $\bar{x}^*(t) = +\infty$. Denote $g_\infty(s):=\lim_{k\downarrow -\infty}g_k(s),s\in [0,T]$. We claim that $ g_\infty(t)= -\infty$. For the sake of contradiction, suppose it is not the case, i.e., $g_\infty(t)>-\infty$. Then, because $g_k(s)$ is decreasing in $s\in [0,T]$ for each $k\le 0$, we have $g_\infty(s)>-\infty,\forall s\in [0,t]$. 
As a result, because $g_\infty$ is monotone on $[0,t]$, it is continuous almost everywhere on $[0,t]$. Recalling that $c_2^*,c_3\in \mathfrak{C}_{\mathrm{pw}}([0, T];\bbR^d)$, that $s_1,s_2\in [\underline{t},t)\cap D$, and that $h(s_1)<h(s_2)$, we can find $\tilde s_1,\tilde s_2 \in (\underline{t},t)\cap D$ with $\tilde s_1<\tilde s_2$ such that $h(\tilde s_1)<h(\tilde s_2)$, that $h$ is continuous at $\tilde s_i,i=1,2$, and that $g_\infty$ is continuous at $\tilde s_i,i=1,2$. Because $g_\infty(\tilde s_1)\ge g_\infty(\tilde s_2)$, there exists $\epsilon_0>0$ such that either $|g_\infty(\tilde s_1)-h(\tilde s_1)|>\epsilon_0$ or $|g_\infty(\tilde s_2)-h(\tilde s_2)|>\epsilon_0$. Without loss of generality, suppose $|g_\infty(\tilde s_1)-h(\tilde s_1)|>\epsilon_0$. Then, if $g_\infty(\tilde s_1)-h(\tilde s_1)>\epsilon_0$, by the continuity of $g_\infty$, $c_3$, and $h$ at $\tilde s_1$, we can find $\delta\in (0,t-\tilde s_1)$ such that $g_\infty(z)-h(z)>\epsilon_0$ and $\|c_3(z)\|\ge \frac{1}{2}\|c_3(\tilde s_1)\|>0$ for all $z\in [\tilde s_1,\tilde s_1+\delta]$, where $\|c_3(\tilde s_1)\|>0$ because $\tilde s_1\in D$. Because $g_k$ is decreasing in $k$, we conclude
\begin{align*}
  g_k(z)-h(z)\ge g_\infty(z)-h(z)>\epsilon_0,\quad \forall z\in [\tilde s_1,\tilde s_1+\delta],\; k<0.
\end{align*}
As a result, for any $k<0$,
\begin{align*}
  &\int_{\underline{t}}^t\|c_3(z)\| | g_k(z)-h(z)|dz\ge \int_{\tilde s_1}^{\tilde s_1+\delta}\|c_3(z)\| | g_k(z)-h(z)|dz\ge \frac{1}{2}\|c_3(\tilde s_1)\|\epsilon_0\delta.
\end{align*}
We then conclude from \eqref{IntegralEquation2} by sending $k$ to $-\infty$ therein that $g_\infty(t) = -\infty$, which contradicts the preassumption that $g_\infty(t)>-\infty$. Thus, we must have $\underline{x}^*(t)=g_\infty(t)=-\infty$. Combining with $\bar{x}^*(t) = +\infty$, we conclude that $\bbX_t = \bbX_t^*=\bbR$.

If $g_\infty(\tilde s_1)-h(\tilde s_1)<-\epsilon_0$, by the continuity of $g_\infty$, $c_3$, and $h$ at $\tilde s_1$, there exists $\delta\in (0,\tilde s_1-\underline{t})$ such that $|g_\infty(z)-g_\infty(\tilde s_1)|+|h(z)-h(\tilde s_1)|<\frac{1}{3}\epsilon_0$ and $\|c_3(z)\|\ge \frac{1}{2}\|c_3(\tilde s_1)\|>0$ for all $z\in [\tilde s_1-\delta,\tilde s_1]$. Moreover, there exists $\bar K>0$ such that $|g_k(\tilde s_1-\delta)-g_\infty(\tilde s_1-\delta)|<\frac{1}{3}\epsilon_0$ for all $k\leq -\bar K$. As a result, because $g_k(s)$ is decreasing in $s\in [0,T]$, we have, for each $z\in [\tilde s_1-\delta, \tilde s_1]$ and $k\leq -\bar K$, that
\begin{align*}
  &g_k(z)-h(z)\le g_k(\tilde s_1-\delta)-h(z)\le g_\infty(\tilde s_1-\delta)+\frac{1}{3}\epsilon_0 -h(z)\\
  & = g_\infty(\tilde s_1-\delta) - g_\infty(\tilde s_1) - \big(h(z)-h(\tilde s_1)) + \frac{1}{3}\epsilon_0 +g_\infty(\tilde s_1) - h(\tilde s_1)<-\frac{1}{3}\epsilon_0.
\end{align*}
As a result,
\begin{align*}
  &\int_{\underline{t}}^t\|c_3(z)\| | g_k(z)-h(z)|dz\ge \int_{\tilde s_1-\delta}^{\tilde s_1}\|c_3(z)\| | g_k(z)-h(z)|dz\ge \frac{1}{6}\|c_3(\tilde s_1)\|\epsilon_0\delta.
\end{align*}
We then conclude from \eqref{IntegralEquation2} by sending $k$ to $-\infty$ therein that $g_\infty(t) = -\infty$, which contradicts the preassumption that $g_\infty(t)>-\infty$. Thus, we must have $\underline{x}^*(t)=g_\infty(t)=-\infty$.

Next, we consider the case in which $h(s)\le 0,\forall s \in [\underline{t},t)$ and $h(s)$ is decreasing in $s\in [\underline{t},t)\cap D$. For each $k\in \bbR$, consider $\tilde g_k$ defined by the following equation:
\begin{align*}
\tilde g_k(s)=\int_{0}^s k\|c_3(z)\| \big(\tilde g_k(z)-h(z)\big)   dz,\quad s\in [0, t].
\end{align*}
Recalling that $c_2^*(z) = 0$ and thus $h(z)=0$ for all $z\in [0,\underline{t})$, we have $\tilde g_k(s)=0,s\in [0,\underline{t}]$ and
\begin{align}\label{eq:tildegk}
  \tilde g_k(s) = -k\int_{\underline{t}}^s e^{k\int_\tau^s\|c_3(z)\|dz   }\|c_3(\tau)\| h(\tau)d\tau,\quad s\in [\underline{t},t].
\end{align}
Fix $k<0$ and consider any $s\in [\underline{t},t)\cap D$. Because $h$ is decreasing on $[\underline{t},t)\cap D$, we have $h(\tau)\ge h(s)$ and thus $\|c_3(\tau)\| h(\tau)\ge \|c_3(\tau)\| h(s)$ for all $\tau\in [\underline{t},s]\cap D$ and $\|c_3(\tau)\| h(\tau)=0\ge \|c_3(\tau)\| h(s)$ for all $\tau\in [\underline{t},s]$ with $c_3(\tau)=0$. As a result, we conclude from \eqref{eq:tildegk} that
\begin{align*}
  \tilde g_k(s)\ge -k\int_{\underline{t}}^s e^{k\int_\tau^s\|c_3(z)\|dz   }\|c_3(\tau)\| h(s)d\tau = h(s) \left(1-e^{k\int_{\underline{t}}^s\|c_3(z)\|dz   } \right) \geq h(s),
\end{align*}
which implies that $\|c_3(s)\|\tilde g_k(s)\ge \|c_3(s)\|\tilde h(s)$. 
As a result, for any $s\in [\underline{t},t)$,
\begin{align*}
  \|c_3(s)\| \big|\tilde g_k(s)-h(s)\big| = \|c_3(s)\| \big(\tilde g_k(s)-h(s)\big),
\end{align*}
so by the uniqueness of the solution to \eqref{IntegralEquation2}, we derive
\begin{align}\label{IntegralEquationSolution}
  g_k(s) = \tilde g_k(s) = -k\int_{\underline{t}}^s e^{k\int_\tau^s\|c_3(z)\|dz   }\|c_3(\tau)\| h(\tau)d\tau,\quad \forall s\in [\underline{t},t].
\end{align}

Fix any $s\in (\underline{t},t]\cap \tilde D$, where $\tilde D:=\{s\in (0,T]: [s-\delta,s)\subseteq D\text{ for some }\delta>0\}$. 
Then, there exists $r_s\in [\underline{t},s)$ such that $c_3(z)\neq 0$ for all $z\in [r_s,s)$. By \eqref{IntegralEquationSolution}, we have
\begin{align*}
  g_k(s) = e^{k\int_{r_s}^s \|c_3(z)\|dz}g_k(r_s)  -k\int_{r_s}^s e^{k\int_\tau^{s}\|c_3(z)\|dz   }\|c_3(\tau)\| h(\tau)d\tau.
\end{align*}
Recall that we already showed that $g_k(z)=\tilde g_k(z) \ge h(z)$ for all $z\in [\underline{t},t)\cap D$ and that $g_k(z)\le 0$ for all $z\in [0,T]$ and $k\le 0$. Also recall that $c_3(z)\neq 0,\forall z\in [r_s,s)$. As a result,
\begin{align*}
 \lim_{k\downarrow -\infty} e^{k\int_{r_s}^s \|c_3(z)\|dz}g_k(r_s) = 0.
\end{align*}
On the other hand, because $h$ is right-continuous and decreasing on $[\underline{t},t)\cap D$, it defines a measure on $(r_s,s)$, so Fubini's theorem yields
\begin{align*}
  -k\int_{r_s}^s e^{k\int_\tau^{s}\|c_3(z)\|dz   }\|c_3(\tau)\| h(\tau)d\tau = h(r_s)\left(1-e^{k\int_{r_s}^s\|c_3(z)\|dz}\right)+ \int_{(r_s,s)} \left(1-e^{k\int_{z}^s\|c_3(\tau)\|d\tau}\right)dh(z).
\end{align*}
Because $\int_{z}^s\|c_3(\tau)\|d\tau>0$ for any $z\in (r_s,s)$, the dominated convergence theorem yields that the limit of the right-hand side of the above equality, as $k$ goes to $-\infty$, is
\begin{align*}
  h(r_s) + \int_{(r_s,s)}  dh(z) = h(s-):=\lim_{z\uparrow s}h(z).
\end{align*}
As a result, $\underline{x}^*(s) = \lim_{k\downarrow -\infty}g_k(s) = h(s-)$.

Now, by the definition of $\tau_t$ and the right-continuity of $c_3$, we have $c_3(z) = 0,\forall z\in [\tau_t,t)$. Because $c_2^*(z)+c_3(z)h(z)=0$ for almost everywhere $z\in [0,\bar t)$ and because $c_2^*$ is right-continuous, we derive $c_2^*(z)=0,\forall z\in [\tau_t,t)$. As a result, $X^*(t) = X^*(\tau_t)$ so $\bbX^*_t = \bbX^*_{\tau_t}$. Thus, $\underline{x}^*(t) = \underline{x}^*(\tau_t) = \lim_{s\uparrow \tau_t}\underline{x}^*(s)$, where the second equality is the case due to the left-continuity of $\underline{x}^*$. By the definition of $\tau_t$ and recalling that $c_2^*$ and $c_3$ are right-continuous, for any $\epsilon\in (0,\tau_t-\underline{t})$, there exists $s\in (\tau_t-\epsilon,\tau_t)$ such that $s\in D$ and that $c_3$ and $h$ are continuous at $s$; in particular, $s\in (\underline{t},t)\cap \tilde D$ and thus $\underline{x}^*(s) = h(s-)=h(s)$. Also recall that $h$ is decreasing on $[\underline{t},t)\cap D$. Then, we conclude $\lim_{s\uparrow \tau_t}\underline{x}^*(s) = \lim_{[\underline{t},\tau_t)\cap D\ni s\uparrow \tau_t}h(s)$, i.e., $\underline{x}^*(t) = \lim_{[\underline{t},\tau_t)\cap D\ni s\uparrow \tau_t}h(s)$. Recalling that $\mathbb{S}_{X(t)} = \lambda_0(t)+\lambda_1(t) \big(x_0+ \mathbb{S}_{X^*(t)}\big)$ and that $\bbX_t = \mathrm{int}(\mathbb{S}_{X(t)})$, we complete the proof of part (iii-b).


Finally, part (iii-c) can be proved similarly.
 \halmos

\subsection{Proof of Corollary \ref{th:PortDensity}}
For any $t\in [0,T]$, because $\sigma(t)\sigma(t)\tran$ is positive definite, so $c^*_2(t)=c_3(t)=0$ if and only if $\theta_0(t)=\theta_1(t)=0$, so \eqref{eq:tupperstarPS} holds.

Next, straightforward calculation yields that for each $\tau\in [0,t^*)$, $\tilde c_2^*(t)+\xi c_3(t)=0,\forall t\in [\tau,t^*)$ if and only if
\begin{align}
  \theta_0(t)  &= \theta_1(t)\left[\int_t^Tb(s)\tran\theta_0(s) e^{-\int_t^sb(z)\tran \theta_1(z)dz}ds - \xi e^{-\int_t^Tb(s)\tran \theta_1(s)ds}\right]\notag\\
   & = \theta_1(t)\left[\int_t^{t^*}b(s)\tran\theta_0(s) e^{-\int_t^sb(z)\tran \theta_1(z)dz}ds - \xi e^{-\int_t^{t^*}b(s)\tran \theta_1(s)ds}\right],\quad t \in [\tau,t^*).\label{eq:theta0Equation}
\end{align}
Because $\theta_1\in \mathfrak{C}_{\mathrm{pw}}([0,T];\bbR^d)$ and $b$ is bounded, \eqref{eq:theta0Equation} has a unique solution. In addition, it is straightforward to verify that
\begin{align*}
  \theta_0(t) = -\xi\theta_1(t),\quad t\in [\tau,t^*)
\end{align*}
solves \eqref{eq:theta0Equation}. As a result, $\tilde c_2^*(t)+\xi c_3(t)=0,\forall t\in [\tau,t^*)$ if and only if $\theta_0(t)+\xi\theta_1(t)=0,\forall t\in[\tau, t^*)$, so \eqref{eq:tlowerstarPS} holds.

Finally, because $\theta_0(s)=\theta_1(s)=0,\forall s\in [t^*,T)$ and $\theta_0(s)+\xi\theta_1(s)=0,\forall s\in[t_*,t^*)$, straightforward calculation that for any $t\in [t_*,t^*)$, $\tilde \xi(t) = \xi$.\halmos

\subsection{Proof of Corollary \ref{prop:PortSupport}}
We first prove \eqref{eq:htildefunPS}, \eqref{eq:tunderlinePS}, and \eqref{eq:tbarPS}. By \eqref{eq:c2starPS} and recalling that $\sigma(t)\sigma(t)\tran$ is positive definite for all $t\in [0,T]$, we conclude that for any $\tau\in (0,T]$, $c_2^*(t)=0,\forall t\in[0,\tau]$ if and only if
\begin{align}\label{eq:theta0Equation1}
  \theta_0(t) = -\theta_1(t)\left(x_0e^{\int_0^tb(s)\tran \theta_1(s) ds } + \int_0^tb(s)\tran \theta_0(s)e^{\int_s^t b(z)\tran\theta_1(z)dz}ds\right),\quad t\in [0,\tau].
\end{align}
Because $\theta_1\in \mathfrak{C}_{\mathrm{pw}}([0,T];\bbR^d)$, the above equation of $\theta_0$ has a unique solution. Moreover, it is straightforward to verify that $\theta_0(t) = -x_0\theta_1(t),t\in[0,\tau]$ solves \eqref{eq:theta0Equation1}. As a result, $c_2^*(t)=0,\forall t\in[0,\tau]$ if and only if $\theta_0(t) = -x_0\theta_1(t),t\in[0,\tau]$. Consequently, we derive \eqref{eq:tunderlinePS}.

Straightforward calculation yields \eqref{eq:hfunPS} and
\begin{align*}
  c_2^*(t) + c_3(t)h(t) = \sigma(t)\tran \big(\theta_0(t) + \theta_1(t)\tilde h(t)\big)e^{-\int_0^tb(s)\tran \theta_1(s)ds},\quad t\in [0,T].
\end{align*}
Because $\sigma(t)\sigma(t)\tran$ is positive definite for all $t\in [0,T]$, we immediately derive \eqref{eq:tbarPS}.

Next, because $\theta_0(t)+x_0\theta_1(t)=0,\forall t\in [0,\underline{t})$, we immediately derive from \eqref{eq:lambdaPS} that $\lambda_0(t) + \lambda_1(t)x_0 = x_0,\forall t\in [0,\underline{t}]$ and from \eqref{eq:xstarPS} that $x^*(t)=x_0,\forall t\in [0,\underline{t}]$. It follows from Theorem \ref{th:GenralSupport}-(i) that $\bbX_t = \{x_0\}$ for all $t\in [0,\underline{t}]$.

Next, we prove that $\mathbb{S}_{X(t)}$ is increasing in $t\in[0,T]$. Lemma \ref{support set multiD} shows that $\mathbb{S}_{X(t)}$ is a closed interval with the lower end $\underline{x}(t)\in[-\infty,+\infty)$ and the upper end $\bar{x}(t)\in(-\infty,+\infty]$. Moreover,
  \begin{align*}
  \underline{x}(t) = \inf_{w\in \mathfrak{H}([0,T];\bbR^d)}f_{w}(t),\quad \bar{x}(t) = \sup_{w\in \mathfrak{H}([0,T];\bbR^d)}f_{w}(t),
\end{align*}
where $f_w$ is given by \eqref{x w multiD} with
\begin{align}\label{eq:SpecializatinPortSel}
  \beta_0 = b\tran \theta_0,\quad \beta_1 = b\tran \theta_1,\quad (\beta_{2,1},\dots, \beta_{2,d})  =\theta_0\tran \sigma, \quad (\beta_{3,1},\dots, \beta_{3,d}) = \theta_1\tran \sigma.
\end{align}
Then, straightforward calculation yields
  \begin{align}\label{fw portfolio}
f_w(t)&=x_0 e^{\int_0^t  \theta_1(s)\tran \big(b(s)-\frac{1}{2}\sigma(s)\sigma(s)\tran \theta_1(s) + \sigma(s)w'(s)\big)ds  } + \int_0^t e^{\int_s^t  \theta_1(\tau)\tran \big(b(\tau)-\frac{1}{2}\sigma(\tau)\sigma(\tau)\tran \theta_1(\tau) + \sigma(\tau)w'(\tau)\big)d\tau  }\notag\\
&\qquad \times \theta_0(s)\tran \big(b(s)-\frac{1}{2}\sigma(s)\sigma(s)\tran \theta_1(s) + \sigma(s)w'(s)\big) ds, \quad t\in [0, T].
\end{align}
For any $s\in[0,T)$, because $\sigma(s)\sigma(s)\tran$ is invertible, for $u(s):=\sigma(s)\tran \left( \frac{1}{2}\theta_1(s)-(\sigma(s)\sigma(s)\tran )^{-1}b(s)  \right)$, we have
\begin{align}\label{no arbitrage}
  b(s)-\frac{1}{2}\sigma(s)\sigma(s)\tran \theta_1(s) + \sigma(s)u(s)=0.
\end{align}
Then, it is straightforward to see that for any $t<s$,
\begin{align*}
  \inf_{w\in \mathfrak{H}([0,T];\bbR^d)}f_{w}(s)\le \inf_{w\in \mathfrak{H}([0,T];\bbR^d)}f_{w}(t),\quad \sup_{w\in \mathfrak{H}([0,T];\bbR^d)}f_{w}(s)\ge \sup_{w\in \mathfrak{H}([0,T];\bbR^d)}f_{w}(t),
\end{align*}
so $\mathbb{S}_{X(t)}\subseteq \mathbb{S}_{X(s)}$.

Next, we already showed that $\bbX_t = \{x_0\}$ for all $t\in [0,\underline{t}]$, so $\bbX_t$ is increasing on $[0,\underline{t}]$. In addition, Theorem \ref{th:GenralSupport} shows that for any $t\in (\underline{t},T]$, $X(t)$ possesses a density function, so $\bbX_t=\mathrm{int}(\mathbb{S}_{X(t)})$. Then, because $\mathbb{S}_{X(t)}$ is increasing in $t\in[0,T]$, we conclude that $\bbX_t$ is increasing in $t\in (\underline{t},T]$. To complete the proof that $\bbX_t$ is increasing in $t\in [0,T]$, we only need to show that $x_0\in \bbX_t$ for any $t\in (\underline{t},T]$.

Consider any $u$ with $u(s)$ solving \eqref{no arbitrage} for any $s\in [0,T]$. For each $K\in \bbR$, construct $w_K\in {\cal U}$ by setting $w_K'(t)=u(t)+K\sigma(t)\tran [\theta_0(t)+x_0 \theta_1(t)], t\in [0, T]$. Then, by \eqref{fw portfolio} and \eqref{no arbitrage}, we derive
  \begin{align*}
f_{w_K}(t)&=  x_0+ K\int_0^t e^{K\int_s^t  \theta_1(\tau)\tran \sigma(\tau) \sigma(\tau) \tran [\theta_0(\tau)+x_0 \theta_1(\tau)]   d\tau  }\|\sigma(s) \tran [\theta_0(s)+x_0 \theta_1(s)]  \|^2   ds,\; t\in [0, T].
\end{align*}
By the definition of $\underline{t}$ and recalling that $\sigma(s)\sigma(s)\tran$ is positive definite for any $s\in[0,T]$, we derive that $f_{w_K}(t)>x_0$ for any $t\in (\underline{t}, T]$ and $K>0$ and $f_{w_K}(t)<x_0$ for any $t\in (\underline{t}, T]$ and $K<0$. Therefore, $\underline{x}(t)<x_0<\bar{x}(t)$, i.e., $x_0$ is in the interior of $\mathbb{S}_{X(t)}$ and thus in $\bbX_t$.

  Finally, assume $\underline{t}<\bar t$ and fix $t\in (\underline{t},\bar t]$ such that $\theta_1(s)\neq 0,s\in (\underline{t},t)$. 
  By the definition of $\bar t$ and recalling \eqref{eq:tbarPS}, we derive $\theta_0(s)+\theta_1(s)\tilde h(s)=0$ for almost everywhere $s\in [\underline{t},t)$ and thus for all $s\in (\underline{t},t)$ because $\theta_0$, $\theta_1$, and $\tilde h$ are right-continuous on $(\underline{t},t)$. Note from \eqref{eq:hfunPS} that $h$ is of finite variation on $(\underline{t},t)$ if and only if $\tilde h$ is of finite variation on $(\underline{t},t)$. Thus, we can differentiate both sides of \eqref{eq:hfunPS} provided that one of $h$ and $\tilde h$ is of finite variation, and the differentiation yields
  \begin{align}
    dh(s) &= e^{-\int_0^sb(z)\tran \theta_1(z)dz}d\tilde h(s) - b(s)\tran \theta_1(s) e^{-\int_0^sb(z)\tran \theta_1(z)dz}\tilde h(s)ds - b(s)\tran \theta_0(s)e^{-\int_0^sb(z)\tran \theta_1(z)dz}ds\notag \\
    & = e^{-\int_0^sb(z)\tran \theta_1(z)dz}d\tilde h(s),\quad s\in (\underline{t},t),\label{eq:htildehrelationPS}
  \end{align}
  where the equality is the case because $\theta_0(s)+\theta_1(s)\tilde h(s)=0,s\in (\underline{t},t)$. 
 Then, we conclude that $h$ is decreasing (increasing, respectively) on $(\underline{t},t)$ if and only if $\tilde h$ is decreasing (increasing, respectively) on $(\underline{t},t)$.

  Now, suppose $\tilde h(s)\le x^*(s),\forall s\in (\underline{t},t)$. Then $h(s)\le 0,\forall s\in (\underline{t},t)$ as $\theta_1(s)\neq 0,s\in (\underline{t},t)$. Because $ h$ is right-continuous at $\underline{t}$ if $\theta_1(\underline{t})\neq 0$ and $h(\underline{t}) =0$ if $\theta_1(\underline{t})=0$, we conclude that $ h(s)\le 0,\forall s\in [\underline{t},t)$. Recall that $h$ is decreasing on $(\underline{t},t)$ if and only if $\tilde h$ is decreasing on $(\underline{t},t)$.
  Then, by Theorem \ref{th:GenralSupport}-(iii)-(b), we immediately conclude that if $\tilde h$ is not decreasing on $(\underline{t},t)$, then $h$ is not decreasing on $(\underline{t},t)$, and thus $\bbX_t = \bbR$.
Next, we consider the case in which   $\tilde h$ is   decreasing on $(\underline{t},t)$ and thus $h$ is decreasing on $(\underline{t},t)$.
  Because $ h$ is right-continuous at $\underline{t}$ if $\theta_1(\underline{t})\neq 0$ and $h(\underline{t}) =0$ if $\theta_1(\underline{t})=0$, and because $ h(s)\le 0,\forall s\in [\underline{t},t)$, we conclude that $h$ is decreasing on $[\underline{t},t)$.
  Then, by Theorem \ref{th:GenralSupport}-(iii)-(b), we immediately conclude that  $\bbX_t = (\lambda_0(t)+\lambda_1(t)(x_0+h(t-)),+\infty)$, where $h(t-):=\lim_{s\uparrow t}h(s)$.
   By \eqref{eq:lambdaPS} and \eqref{eq:hfunPS}, straightforward calculation yields $\lambda_0(t)+\lambda_1(t)(x_0+h(t-))=\tilde h(t-)$. Thus, $\bbX_t = (\tilde h(t-),+\infty)$.

  Moreover, if $\tilde h$ is decreasing on $(\underline{t},t)$, we claim that $\tilde h(s)<x^*(s),\forall s\in (\underline{t},t)$. For the sake of contradiction, suppose $\tilde h(s_0)=x^*(s_0)$ for some $s_0\in (\underline{t},t)$. Then, by \eqref{eq:hfunPS}, we have $h(s_0)=0$. Recall we have shown that  $ h(s)\le 0,\forall s\in [\underline{t},t)$, and  $h$ is decreasing on $(\underline{t},t)$ because $\tilde h$ is   decreasing on $(\underline{t},t)$.
  As a result, $h(s)=0,\forall s\in (\underline{t},s_0]$ and, consequently, $\tilde h(s) =x^*(s),\forall s \in (\underline{t},s_0]$.
    It follows from \eqref{eq:htildehrelationPS} that $dx^*(s) = d\tilde h(s) = 0,\forall s\in (\underline{t},s_0)$. Because $x^*$ is continuous on $[0,T]$ and $x^*(s) =x_0,\forall s\in [0,\underline{t}]$, we conclude that $x^*(s) = x_0,\forall s \in [0,s_0]$. As a result, $\tilde h(s) =x^*(s)= x_0,\forall s\in (\underline{t},s_0)$. By \eqref{eq:tbarPS}, we have $\theta_0(s) + \tilde h(s) \theta_1(s)=0$ for almost everywhere $s\in[0,\bar t]$. Together with the right continuity of $\theta_0$, $\theta_1$, and $\tilde h$ on $(\underline{t},t)$, we derive $\theta_0(s) + x_0\theta_1(s)=\theta_0(s) + \tilde h(s) \theta_1(s)=0$ for all $s\in (\underline{t},s_0]$. Because $\theta_0$ and $\theta_1$ are right-continuous, we derive $\theta_0(s) + x_0\theta_1(s)=0$ for all $s\in [\underline{t},s_0]$. This contradicts the definition of $\underline{t}$.


  Finally, the case in which $\tilde h(s)\ge  x^*(s),\forall s\in (\underline{t},t)$ can be treated similarly.
\halmos

\appendix

\bibliography{LongTitles,BibFile}

\end{document}